%% file: main.tex
\documentclass[a4paper,cleveref,autoref,thm-restate]{lipics-v2021}
\usepackage{cite}
\usepackage{amsmath,amsthm,amssymb,amsfonts}
\usepackage[linesnumbered,ruled,vlined]{algorithm2e}
\usepackage{threeparttable}
\usepackage{tabularx}
\usepackage{graphicx}
\usepackage{subcaption}
\usepackage{textcomp}
\usepackage{multirow}
\usepackage[dvipsnames,table,xcdraw]{xcolor}
\usepackage{tkz-graph}
\usepackage{makecell}
\usepackage{bbm}
\usepackage{appendix}
\usepackage{stfloats}
\usepackage{float}
\usepackage{xspace}
\usepackage{comment}
\usepackage[normalem]{ulem}
\usepackage{tikz}
\usepackage{todonotes}

\SetCommentSty{mycommfont}

\newtheorem{prune}{Pruning Rule}

\title{Multi-Objective Memory Bandwidth Regulation and Cache Partitioning for Multicore Real-Time Systems (with Supplemental Material)}

\titlerunning{Multi-Objective Memory Bandwidth Regulation and Cache Partitioning}

\author{Binqi Sun}{Technical University of Munich, Germany}{binqi.sun@tum.de}{https://orcid.org/0000-0002-9764-6259}{}
\author{Zhihang Wei}{Technical University of Munich, Germany}{zhihang.wei@tum.de}{https://orcid.org/0000-0003-3755-3509}{}
\author{Andrea Bastoni}{Technical University of Munich, Germany}{andrea.bastoni@tum.de}{https://orcid.org/0000-0001-8256-6160}{}
\author{Debayan Roy}{Technical University of Munich, Germany}{debayan.roy.tum@gmail.com}{}{}
\author{Mirco Theile}{Technical University of Munich, Germany}{mirco.theile@tum.de}{https://orcid.org/0000-0003-1574-8858}{}
\author{Tomasz Kloda}{LAAS-CNRS, Insa de Toulouse, France}{tkloda@laas.fr}{https://orcid.org/0000-0003-0822-4976}{}
\author{Rodolfo Pellizzoni}{University of Waterloo, Canada}{rpellizz@uwaterloo.ca}{https://orcid.org/0000-0002-7331-804X}{}
\author{Marco Caccamo}{Technical University of Munich, Germany}{mcaccamo@tum.de}{https://orcid.org/0000-0003-2328-044X}{}



\nolinenumbers 

\EventEditors{Renato Mancuso}
\EventNoEds{1}
\EventLongTitle{37th Euromicro Conference on Real-Time Systems (ECRTS 2025)}
\EventShortTitle{ECRTS 2025}
\EventAcronym{ECRTS}
\EventYear{2025}
\EventDate{July 8--11, 2025}
\EventLocation{Brussels, Belgium}
\EventLogo{}
\SeriesVolume{335}
\ArticleNo{7}

\relatedversion{}

\Copyright{Binqi Sun, Zhihang Wei, Andrea Bastoni, Debayan Roy, Mirco Theile, Tomasz Kloda, Rodolfo Pellizzoni, and Marco Caccamo}

\authorrunning{B. Sun et al.}

\acknowledgements{Marco Caccamo was supported by an Alexander von Humboldt Professorship endowed by the German Federal Ministry of Education and Research.}

\input{macros}

\begin{document}

\maketitle

\begin{abstract}
    Memory bandwidth regulation and cache partitioning are widely used techniques for achieving predictable timing in real-time computing systems. Combined with partitioned scheduling, these methods require careful co-allocation of tasks and resources to cores, as task execution times strongly depend on available allocated resources. 
    To address this challenge, this paper presents a 0-1 linear program for task-resource co-allocation, along with a multi-objective heuristic designed to minimize resource usage while guaranteeing schedulability under a preemptive EDF scheduling policy. 
    Our heuristic employs a multi-layer framework, where an outer layer explores resource allocations using Pareto-pruned search, and an inner layer optimizes task allocation by solving a knapsack problem using dynamic programming. 
    To evaluate the performance of the proposed optimization algorithm, we profile real-world benchmarks on an embedded AMD UltraScale+ ZCU102 platform, with fine-grained resource partitioning enabled by the Jailhouse hypervisor, leveraging cache set partitioning and MemGuard for memory bandwidth regulation. 
    Experiments based on the benchmarking results show that the proposed 0-1 linear program outperforms existing mixed-integer programs by finding more optimal solutions within the same time limit. Moreover, the proposed multi-objective multi-layer heuristic performs consistently better than the state-of-the-art multi-resource-task co-allocation algorithm in terms of schedulability, resource usage, number of non-dominated solutions, and computational efficiency. 
\end{abstract}

\input{sec_introduction}

\input{sec_literature}

\input{sec_system_model}

\input{sec_method}

\input{sec_tools}

\input{sec_experiments}

\input{sec_conclusion}

\bibliographystyle{plainurl}
\bibliography{reference}

\input{supplemental_material}

\end{document}

%% file: macros.tex
\usepackage[]{todonotes}
\definecolorseries{test}{rgb}{grad}[rgb]{.95,.55,.55}{11,11,17}
\resetcolorseries[10]{test}
\newcommand{\addtodoeditor}[1]{%
    \colorlet{#1}{test!!+!50}
    \expandafter\newcommand\csname#1\endcsname [1]{%
        \todo[color=#1,size=\tiny]{\sffamily\textbf{\uppercase{#1}:}
    ##1}\xspace%
    }
    \expandafter\newcommand\csname#1i\endcsname [1]{%
        \todo[inline, color=#1]{\sffamily\textbf{\uppercase{#1}:} ##1}\xspace%
    }
}

\addtodoeditor{bs}
\addtodoeditor{dr}
\addtodoeditor{ab}
\addtodoeditor{tk}
\addtodoeditor{rp}
\addtodoeditor{mt}

\newcommand{\eg}{{\it e.g.}\xspace}
\newcommand{\ie}{{\it i.e.}\xspace}

\newcommand{\etal}{{\it et al.}\xspace}

\let\oldnl\nl
\newcommand{\nonl}{\renewcommand{\nl}{\let\nl\oldnl}}

\newcommand{\sdvbshort}{\mbox{\emph{SD-VBS}}\xspace}
\newcommand{\sdvblong}{\emph{San Diego Vision Benchmark Suite}\xspace}
\newcommand{\parsec}{\mbox{\emph{PARSEC}}\xspace}

%% file: sec_introduction.tex
\section{Introduction}
\label{sec:intro}

The advent of multiprocessor systems-on-chip (MPSoC) has transformed the landscape of embedded real-time computing by enabling platforms that combine high performance with energy efficiency. Modern MPSoCs, such as the AMD Ultrascale+ ZCU102~\cite{zcu-102} and NVIDIA Orin~\cite{nvidia-agx-orin}, integrate multiple processing cores and specialized accelerators within a single system, sharing a unified memory hierarchy. 
These platforms are essential in addressing the increasing demand for complex real-time applications. However, fully utilizing their potential presents significant challenges due to cross-core interference arising from the shared nature of critical resources such as memory interconnects, memory controllers, and caches, which complicates the analysis required for real-time guarantees~\cite{guan2009cache,dong2018shared,maiza2019survey,davis2022framework,lin2023lag,lesage2023reducing}. 

To address these challenges, both hardware- and software-based partitioning techniques have been developed to manage the allocation of shared resources.
Hardware solutions, such as Intel RDT~\cite{intel-rdt} and Arm MPAM~\cite{arm-mpam}, offer 
monitoring and partitioning capabilities for hardware resources like caches and interconnects to ensure quality of service (QoS) levels.
Meanwhile, software-based techniques, including memory bandwidth regulation~\cite{yun2013memguard, MemPol} and cache partitioning~\cite{MDBCCP:13, Kim16:EMSOFT, KWCFAS:17}, have become popular because they offer flexibility across multiple platforms and can be integrated at the hypervisor level. This integration ensures transparency to both operating systems and applications.
By controlling the allocation of shared resources, these partitioning techniques improve isolation and predictability, enabling tighter bounds on worst-case execution times (WCETs) and response time analysis for real-time tasks \cite{lv2016survey,XPCLLLL:19}. 
However, they also introduce new complexities in task and resource co-allocation. 
Task execution times can be affected significantly by the availability of shared resources, necessitating careful coordination to balance resource usage and task allocation for high system performance and real-time schedulability. 

This work addresses the problem of task and resource co-allocation in the context of partitioned scheduling, which was shown to be effective in practice (\eg,~\cite{DBLP:conf/rtss/BrandenburgG16}).
Specifically, we focus on systems employing a preemptive Earliest Deadline First (EDF) scheduling policy. The objective is to minimize memory bandwidth and cache resource usage while ensuring real-time schedulability. We formulate this problem as a multi-objective zero-one linear program (0-1 ILP). Although the ILP can be solved using standard mathematical programming solvers to derive optimal solutions for each single objective separately, it suffers from an exponential time complexity and cannot generate multiple Pareto-optimal solutions within a reasonable time. To address these challenges, we propose a multi-objective resource-task co-allocation heuristic that employs a multi-layer optimization framework. The outer layer explores resource allocations using Pareto-pruned search, and the inner layer optimizes task allocation by solving a knapsack problem with dynamic programming. 

To evaluate our approach, we derive \emph{slowdown profiles} for real-world benchmarks from the \sdvblong (\sdvbshort)~\cite{venkata2009sd, rt-bench} and \parsec~\cite{parsec} suites, executed on a real embedded MPSoC platform (AMD UltraScale+ ZCU102). These profiles quantify execution time slowdowns as a function of allocated  memory bandwidth and cache partitions (compared to the configuration with full resource availability).
We conducted extensive benchmarking on the ZCU102, which features four Cortex-A53 cores sharing a 1~MB last-level cache. Specifically, we configured fine-grained resource partitioning schemes using the Jailhouse hypervisor~\cite{jailhouse, minerva-jailhouse}, leveraging cache set partitioning (coloring, \eg,~\cite{KSMCV:19}) and MemGuard~\cite{yun2013memguard} for bandwidth regulation. 
Our experiments based on these real-world benchmarks demonstrate that the proposed 0-1 ILPs are solved more efficiently than an existing model based on mixed-integer programming (MIP)~\cite{XPCLLLL:19}. Moreover, our proposed heuristic consistently outperforms the state-of-the-art algorithm in terms of schedulability, resource usage, number of non-dominated solutions, and solving efficiency.

In summary, we make the following contributions:
\begin{itemize}
    \item We develop a 0-1 ILP for the resource-task co-allocation problem, outperforming the existing mixed-integer program developed in~\cite{XPCLLLL:19} that includes non-binary decision variables and non-linear constraints (Section~\ref{sec:sys_model}). 
    \item We propose a multi-objective multi-layer optimization heuristic (MMO) with effective Pareto-based pruning strategies and knapsack-based task allocation to minimize memory bandwidth and cache usage simultaneously while ensuring schedulability (Section~\ref{sec:method}). 
    \item We profile real-world benchmarks from two test suites (\sdvbshort~\cite{venkata2009sd} and \parsec~\cite{parsec}) on an embedded MPSoC hardware to evaluate the impact of the memory bandwidth and cache allocation on task execution slowdowns (Section~\ref{sec:tools}). 
    \item We demonstrate the effectiveness and efficiency of the proposed methods by comparing them with state-of-the-art methods through extensive experimental evaluations using real-world benchmarks. Specifically, MMO identifies up to 62.67\% and 50.34\% more schedulable task sets than the state-of-the-art algorithm on the \sdvbshort and \parsec benchmarks, respectively, with significant memory bandwidth and cache usage savings and computation time reduction (Section~\ref{sec:exp}). 
\end{itemize}

%% file: sec_literature.tex
\section{Related Work}
\label{sec:literature}

In this section, we give an overview of the related works on task and resource allocation strategies for real-time systems and discuss the differences between our proposed algorithm and the state-of-the-art resource-task co-allocation methods. 

\subsection{Task Allocation}
Mapping tasks statically to individual processors is widely used in industry practice due to its low scheduling overhead~\cite{DBLP:conf/rtss/BrandenburgG16}. 
However, since the task allocation problem is NP-hard in the strong sense~\cite{ekberg2021partitioned}, many approximation methods have been developed for both 
preemptive~\cite{Burchard:1995,Dhall:1978,Baruah:2005,Lopez:2000} and non-preemptive~\cite{Fisher:2006,Senoussaoui:2020} scheduling policies. 
These methods have also been extended to support parallel task scheduling, including directed acyclic graphs (DAGs)~\cite{fonseca2016response,casini2018partitioned,Zahaf:2020} and gang tasks~\cite{Ueter:2021,sun2024strict,sun2024partitioned}, as well as to take into account inter-task interference~\cite{Zahaf:2021}.
On the other hand, exact approaches to the partitioning problem use optimization techniques such as mixed-integer linear programming (MILP)~\cite{abeni2022partitioning,Mo:2023}. While MILP formulations can provide exact solutions, their scalability remains a challenge, particularly in systems with a large number of tasks or processors. 
A detailed discussion on the precise complexity classes of a list of real-time task allocation problems can be found in~\cite{ekberg2021partitioned}. 

\subsection{Resource Allocation}
Cache and memory bandwidth are two critical resources to be partitioned for achieving timing predictability in real-time systems. 
A widely adopted \emph{software-based} approach to cache partitioning is \emph{cache coloring}, which has been implemented at both the operating system (OS)~\cite{MDBCCP:13, Kim16:EMSOFT, KWCFAS:17} and hypervisor levels~\cite{KSMCV:19, xilinx-xen-cache-color}. In this paper, we rely on a cache-coloring implementation available in the Jailhouse hypervisor~\cite{minerva-jailhouse}. 
Alternatively, caches can also be partitioned via hardware modifications (\eg,~\cite{Survey-Way-Part}) or by exploiting hardware support such as the Arm DSU~\cite{arm-dynamiciq}, which is only available on very recent embedded Arm platforms and notably not yet supported on our Ultrascale+.
Similarly to caches, memory bandwidth partitions can be assigned in software leveraging hardware features such as Performance Monitoring Units (PMUs). For example, MemGuard~\cite{yun2013memguard} and MemPol~\cite{MemPol} propose a per-core memory bandwidth partitioning using PMU-based counters.
Hardware modifications to generally improve the predictability of memory accesses have also been proposed (\eg,~\cite{hassan2019reduced, BRU:20}).
Intel RDT~\cite{intel-rdt} supports partitioning of both caches and memory bandwidth and has been used in \eg,~\cite{XPCLLLL:19}. Nonetheless, real-time characteristics of Intel RDT have been found to be not always effective~\cite{SBMYK:22}.
Arm MPAM~\cite{arm-mpam} is a recent specification with partitioning capabilities similar to Intel RDT, but to date, no available implementations for COTS platforms exist.

Building on these cache and memory bandwidth partitioning methods, various allocation strategies have been developed to effectively dedicate resources to real-time tasks and improve schedulability. For caches, approaches such as branch-and-bound~\cite{Altmeyer:2014,Altmeyer:2016}, genetic algorithms~\cite{Bui:2008,Meroni:2023}, and guided-local search~\cite{sun2023minimizing,sun2024minimizing} have been proposed to optimize how cache partitions are assigned to real-time workloads.
Similar efforts exist for memory bandwidth allocation. Aghilinasab~\etal~\cite{aghilinasab2020dynamic} present a dynamic scheme that monitors and reallocates memory bandwidth between real-time and best-effort tasks, adapting to runtime variations. Park~\etal~\cite{park2019copart} further propose a coordinated approach for LLC and memory bandwidth partitioning, targeting workload fairness rather than hard timing guarantees.

\subsection{Task and Cache Co-Allocation}
Beyond independent allocation strategies for real-time tasks and resources, the co-allocation of tasks and cache partitions has been explored to further enhance real-time schedulability.
Under preemptive EDF scheduling, Chisholm~\etal~\cite{CWKA:15} introduce MC$^2$, a linear programming-based optimization framework for mixed-criticality multicore real-time systems. 
Kim and Rajkumar~\cite{Kim16:EMSOFT} develop a cache management scheme for cache-to-task allocation and later proposed a cache-aware task allocation algorithm tailored for virtual machine design.
The task and cache allocation under non-preemptive scheduling introduces additional challenges due to blocking effects, making task utilization an insufficient sole indicator for schedulability. Berna and Puaut~\cite{Berna:2012} propose a period-driven task and cache partitioning algorithm under non-preemptive EDF scheduling, prioritizing task period compatibility as the primary partitioning criterion. Paolieri~\etal~\cite{Paolieri:2011} introduce IA$^3$, an interference-aware allocation algorithm focusing on WCET sensitivity. Sun~\etal~\cite{sun2023co} develop a search-based algorithm leveraging a first-fit heuristic for task allocation and propose two heuristic variants considering task period compatibility and cache sensitivity as the task ordering criteria. 
However, these works do not consider memory bandwidth allocation. Ignoring contention on the shared memory bus can compromise system predictability while assuming uniform memory bandwidth partitioning limits the flexibility needed to optimize schedulability effectively. 

\subsection{Task and Multi-Resource Co-Allocation}
\label{sec:related_work_multi_resource_co_alloc}
Recent studies~\cite{XPCLLLL:19,Meng:2019,Nie:2022,gifford2021dna} have also explored the problem of task and multi-resource co-allocation. 
Meng \etal~\cite{Meng:2019} propose a multi-resource allocation framework for real-time multicore virtualization, incorporating techniques to mitigate abstraction overhead. 
Nie \etal~\cite{Nie:2022} investigate federated scheduling for parallel tasks, where each task is assigned a set of cores along with cache and memory bandwidth partitions while ensuring schedulability conditions. 
Gifford \etal~\cite{gifford2021dna} employ a worst-fit bin-packing approach to allocate soft real-time tasks to cores, dynamically adjusting cache and memory bandwidth partitions based on task deadlines.

The most \emph{closely related work} to ours is \cite{XPCLLLL:19}, which also addresses the co-allocation of tasks and resources, including memory bandwidth and cache partitions. Their algorithm, CaM, represents the current state-of-the-art and outperforms other previous methods. It employs a k-means algorithm to classify tasks into different clusters, followed by a task-cluster-to-core allocation heuristic based on first-fit bin-packing and a resource-to-core allocation heuristic based on resource utility. Our method, MMO, improves upon CaM in three key ways:
\begin{itemize}
    \item The multi-layer search strategy of MMO provides a more thorough exploration of task and resource co-allocation possibilities. 
    \item The inner-layer task allocation in MMO is formulated as a knapsack problem, enabling the use of dynamic programming to achieve optimal task assignments that maximize core utilization. 
    \item MMO is a multi-objective heuristic that simultaneously optimizes the usage of multiple resources, yielding multiple non-dominated solutions, while CaM can only produce a single solution. 
\end{itemize}
In Section~\ref{sec:exp}, we present extensive experimental evaluations to compare the performance of MMO with CaM in terms of schedulability, resource usage, and computation efficiency.

%% file: sec_system_model.tex
\section{System Model and Mathematical Formulation}
\label{sec:sys_model}

\subsection{System Model}

We consider a set of $N$ real-time tasks $\tau$ running on $M$ cores that share a last-level cache (LLC) and memory bus. 
The memory bandwidth and the LLC are divided into $B$ and $K$ uniform partitions, respectively.
Each task $\tau_i \in \tau$ generates an infinite sequence of sporadic job instances, with a minimum inter-arrival time of $T_i$, and each job must complete its execution before the next instance arrives (\ie, has an implicit deadline). 
Tasks are scheduled by a partitioned scheduling policy, where the task set is divided into $M$ subsets, each statically assigned to a specific core (\ie, jobs of a task remain on their assigned core without migrating at runtime). 
Within each core, tasks are scheduled by a preemptive Earliest Deadline First (EDF) scheduler. 
Partitions of shared resources (\ie, memory bandwidth and LLC) are also statically assigned to each core, providing inter-core isolation to reduce interference. 
The worst-case execution time (WCET) of a task $\tau_i$, given $b$ memory bandwidth and $k$ cache partitions, is characterized by $C_i(b, k)$. 
We assume that the cache-related preemption delays (CRPD) are accounted for in the task WCETs.
The reference WCET, $\widehat{C}_i = C_i(K, B)$, gives the WCET when the task has access to full memory bandwidth and LLC. 
Task utilization is defined as $U(\tau_i, b,k) = C_i(b,k) / T_i$, with $\widehat{U}(\tau_i) = \widehat{C}_i / T_i$ as the reference utilization. 
Under the EDF scheduling policy~\cite{liu1973scheduling}, a subset of tasks executing on the same core, denoted as $\tau'$, is schedulable if and only if the sum of their utilizations does not exceed 1 (\ie, $\sum_{\tau_i\in\tau'}{U(\tau_i,b,k)} \leq 1$), given $b$ memory bandwidth and $k$ cache partitions allocated to the core. 
The slowdown of a task under $b$ memory bandwidth and $k$ cache allocation is defined as the ratio of its WCET compared to its reference WCET, \ie, $C_i(b, k) / \widehat{C}_i$.

The slowdown of an example benchmark, \texttt{mser-qcif} from the \sdvbshort test suite, is visualized in Figure~\ref{fig:bench}. The upper inset of the figure presents the slowdown profile as a heat map for various numbers of $b$ memory bandwidth and $k$ cache partitions, with the numbers at the corners representing the slowdowns under extreme resource allocations. The bottom insets explicitly visualize the slowdown trends when increasing one of the controlled dimensions ($b$ or $k$) while keeping the other constant. 
The details of the profiling experiments and more benchmark profiles are presented in Section~\ref{sec:tools}. 

A summary of the notations used in the paper is given in Table~\ref{tab:notation}.

\input{figures/example_bench}

\input{tables/notions}

\subsection{Mathematical Formulation}
As shown in Figure~\ref{fig:bench}, the task execution times can be significantly influenced by the resources (\ie, memory bandwidth and cache partitions) allocated to the core on which they execute. Therefore, ensuring schedulability on all cores requires a careful co-allocation of tasks and resources such that no task misses its deadline. Moreover, simply meeting deadlines is often insufficient in embedded systems, where resources are inherently limited and frequently shared among multiple applications. Because these constrained platforms must also address cost and power limitations, minimizing resource usage becomes a critical objective.

By reducing the resources reserved for real-time tasks to only what is strictly necessary, system designers can allocate remaining resources to best-effort tasks or other value-added services without jeopardizing the timing guarantees of critical workloads. In some cases, disabling unused resources leads to lower power consumption and reduced heat dissipation, thus extending battery life and enhancing overall system reliability~\cite{armbw, mittal2014survey, david2011memory}. This balance between schedulability and resource efficiency cannot be captured by a feasibility formulation focusing solely on meeting deadlines. 

A multi-objective optimization approach naturally addresses these competing demands by finding configurations that simultaneously guarantee real-time performance and minimize multiple resource usage. In this way, system designers gain a set of Pareto-optimal (non-dominated) solutions--each reflecting a distinct trade-off between the usage of different resources. Such flexibility is highly beneficial in the diverse and evolving landscape of real-time applications, where a single, one-size-fits-all solution often falls short.

Based on these considerations, we formulate a multi-objective optimization problem as follows with notations in Table~\ref{tab:notation}. 

\begin{align}
    &\text{min.} &&\sum\nolimits_{m=1}^{M}\sum\nolimits_{b=1}^{B}{y_{bm}} \cdot b \label{eq:obj_b},\\
    &\text{min.} &&\sum\nolimits_{m=1}^{M}\sum\nolimits_{k=1}^{K}{z_{km}} \cdot k \label{eq:obj_k},\\
    &\text{s.t.} &&\sum\nolimits_{m=1}^{M}{x_{im}} = 1, \quad \forall i \in [1,N],\label{eq:const_task_alloc}\\
    & &&\sum\nolimits_{b=1}^{B}{y_{bm}} = 1, \quad \forall m \in [1,M],\label{eq:const_b_alloc}\\
    & &&\sum\nolimits_{k=1}^{K}{z_{km}} = 1, \quad \forall m \in [1,M],\label{eq:const_k_alloc}\\
    & &&\sum\nolimits_{m=1}^{M}\sum\nolimits_{b=1}^{B}{y_{bm}} \cdot b \leq B,\label{eq:const_b_total}\\
    & &&\sum\nolimits_{m=1}^{M}\sum\nolimits_{k=1}^{K}{z_{km}} \cdot k \leq K,\label{eq:const_k_total}\\
    & &&3 \cdot \alpha_{ibkm} \leq x_{im} + y_{bm} + z_{km}, \quad \forall i \in [1,N], b \in [1,B], k \in [1,K], m \in [1,M],\label{eq:const_alpha_1}\\
    & &&\alpha_{ibkm} \geq x_{im} + y_{bm} + z_{km} - 2, \quad \forall i \in [1,N], b \in [1,B], k \in [1,K], m \in [1,M],\label{eq:const_alpha_2}\\
    & &&\sum\nolimits_{i=1}^{N}\sum\nolimits_{b=1}^{B}\sum\nolimits_{k=1}^{K}{\alpha_{ibkm} \cdot U(\tau_i,b,k)} \leq 1, \quad \forall m \in [1,M],\label{eq:const_util}
\end{align}
where \eqref{eq:obj_b} and \eqref{eq:obj_k} define multi-objective functions to minimize memory bandwidth and cache usage, respectively. Constraints~\eqref{eq:const_task_alloc} ensure each task is assigned to exactly one core, while \eqref{eq:const_b_alloc} and \eqref{eq:const_k_alloc} guarantee each core is assigned a specific number of memory bandwidth and cache partitions, respectively. Constraints~\eqref{eq:const_b_total} and \eqref{eq:const_k_total} enforce the total allocated memory bandwidth and cache partitions do not exceed their capacities. Constraints \eqref{eq:const_alpha_1} and \eqref{eq:const_alpha_2} provide linear equivalents for the conjunction $x_{im} \land y_{bm} \land z_{km}$, where $\alpha_{ibkm}=1$ if task $i$ is assigned to core $m$, which allocates $b$ memory bandwidth and $k$ cache partitions. 
Constraints~\eqref{eq:const_util} ensure that the assigned tasks are schedulable (\ie, satisfy the utilization bound) on each core.

The resulting model is a zero-one linear program (0-1 ILP) since it contains only binary variables and linear constraints. 
In Section~\ref{sec:eval_results}, we will show that our 0-1 ILP model is solved more efficiently than the mixed-integer program (MIP) developed by \cite{XPCLLLL:19}, which includes non-binary decision variables and non-linear constraints. 

\subsection{Optimization Challenges}
The 0-1 ILP can be solved using a standard mathematical programming solver to optimize every single objective, \ie, memory bandwidth~\eqref{eq:obj_b} or cache partitions~\eqref{eq:obj_k}, separately. However, there remain two significant challenges to be addressed: 
\begin{itemize}
    \item \textbf{Scalability.} Even though the 0-1 ILP only contains linear constraints, the solution space grows exponentially with respect to the number of tasks, cores, and resource partitions, making it challenging for solvers to find feasible solutions within a reasonable time. This limits its practicality in real-world system design.
    \item \textbf{Multi-objective optimization.} The formulation requires minimizing the usage of multiple resources simultaneously. 
    While multi-objective optimization problems can be transformed into multiple single-objective problems using a weighted combination of objectives, the resulting solutions do not necessarily correspond to the pre-determined weights, failing to reflect the desired trade-off~\cite{deb2011multi}. 
    A more effective approach is to develop a multi-objective optimization algorithm that generates a set of Pareto-optimal solutions, where no solution dominates another in terms of all objectives. This allows system designers to choose the most suitable allocation from the Pareto set based on their specific preferences and trade-offs. 
\end{itemize}
Therefore, we propose a new heuristic to address these challenges.

%% file: figures/example_bench.tex
\begin{figure}[H]
    \centering
    \begin{subfigure}[]{\textwidth}
        \centering
        \begin{tikzpicture}[]
            \node[] at (1,1) {\includegraphics[width=5cm]{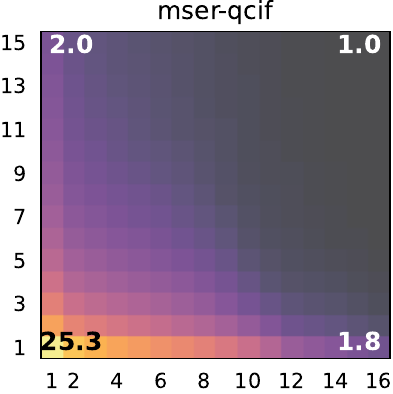}};
            
            \node[] at (4.5,1) {\includegraphics[height=5cm]{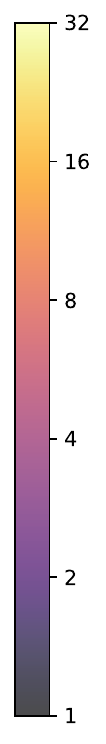}};        
    
            \node[rotate=90, align=center] at (-2, 1) {Number of Memory \\ Bandwidth Partitions};
            \node[] at (1.2, -1.8) {Number of Cache Partitions};
            \node[rotate=90] at (5, 1) {Slowdown (log-scale)};
        \end{tikzpicture}
    \end{subfigure}
    \\
    \begin{subfigure}[]{0.45\textwidth}
        \begin{tikzpicture}[]
            \node[] at (0,1) {\includegraphics[width=7cm]{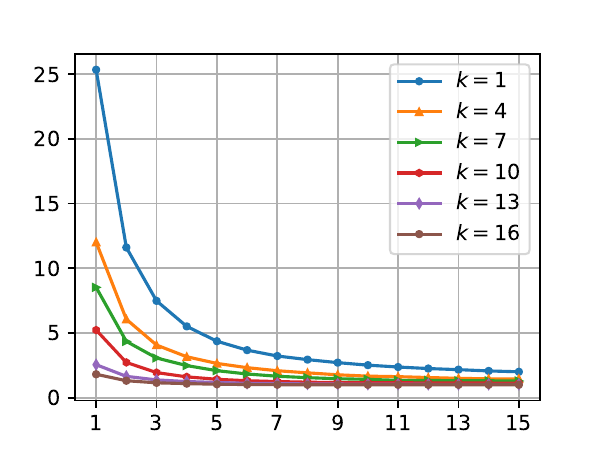}};
    
            \node[rotate=90, align=center] at (-3.5, 1) {Slowdown};
            \node[] at (0, -1.8) {Number of Memory Bandwidth Partitions};
        \end{tikzpicture}
    \end{subfigure}
    \hspace{0.5cm}
    \begin{subfigure}[]{0.45\textwidth}
        \begin{tikzpicture}[]
            \node[] at (0,1) {\includegraphics[width=7cm]{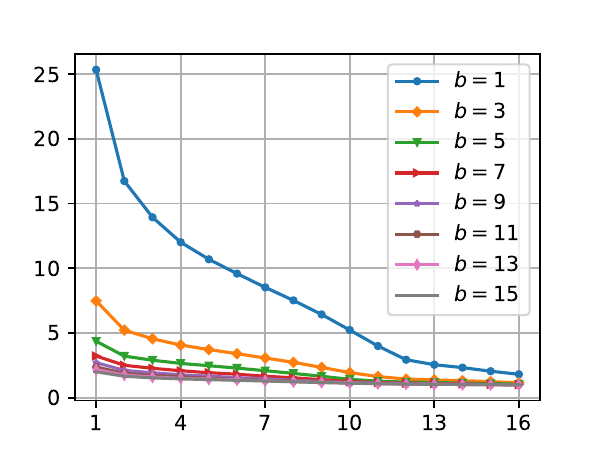}};
    
            \node[rotate=90, align=center] at (-3.5, 1) {Slowdown};
            \node[] at (0, -1.8) {Number of Cache Partitions};
        \end{tikzpicture}
    \end{subfigure}
    
    \caption{Slowdown profile (relative to unrestricted execution) of the \texttt{mser-qcif} benchmark with variable memory bandwidth $b$ and cache partition $k$ allocations. The upper inset presents the slowdown as a heatmap where both dimensions vary; the bottom insets present the resource sensitivity of each dimension separately.} 
    \label{fig:bench}
\end{figure}

%% file: tables/notions.tex
\begin{table}[t]
  \begin{center}
    \caption{Summary of notations.}
    \label{tab:notation}
    \renewcommand{\arraystretch}{1.15}
    \begin{tabular}{p{30mm}p{100mm}}
      \hline
      \textbf{Notations} & \textbf{Description} \\ \hline 
      \multicolumn{2}{l}{System model} \\
      $\tau$ & a set of $N$ tasks $\tau=\{\tau_1,\tau_2,\dots,\tau_{N}\}$; \\ 
      $M, B, K$ & total number of cores, memory bandwidth, and cache partitions; \\ 
      $U(\tau_i, b, k)$ & $\tau_i$'s utilization with $b$ memory bandwidth and $k$ cache partitions;  \\
      $\widehat{U}(\tau_i)$ & $\tau_i$'s reference utilization;  \\ \hline
      \multicolumn{2}{l}{Mathematical formulation} \\
      $x_{im}$ & 1, if $\tau_i$ is assigned to core $m$; 0, otherwise \\
      $y_{bm}$ & 1, if $b$ memory bandwidth partitions are assigned to core $m$; 0, otherwise \\
      $z_{km}$ & 1, if $k$ cache partitions are assigned to core $m$; 0, otherwise \\
      $\alpha_{ibkm}$ & auxiliary variables: $\alpha_{ibkm} = x_{im} \land y_{bm} \land z_{km}$ \\ \hline
      \multicolumn{2}{l}{Optimization heuristic} \\
      $\omega,\omega'$ & a (new) partial solution; \\
      $\omega.\overrightarrow{\tau}, \omega.\overrightarrow{b}, \omega.\overrightarrow{k}$ & $\omega$'s task, memory bandwidth, and cache allocation vectors; \\ 
      $\omega.\overline{\tau}, \omega.\overline{B}, \omega.\overline{K}$ & $\omega$'s remaining tasks, memory bandwidth, and cache partitions; \\ 
      $\Omega, \Omega'$ & a (new) set of partial solutions; \\
      $\Omega^*$ & a set of complete non-dominated solutions; \\ \hline
    \end{tabular}
  \end{center}
  \vspace{-1.5em}
\end{table}

%% file: sec_method.tex
\section{Multi-Objective Multi-Layer Optimization}
\label{sec:method}

This section introduces MMO, a multi-objective, multi-layer optimization heuristic for task-resource co-allocation. 
As discussed in Section~\ref{sec:sys_model}, task and resource allocations are inherently interdependent, necessitating a thorough exploration of the joint design space, which is known to be NP-hard~\cite{XPCLLLL:19}. To address this challenge, we propose a multi-layer design space exploration heuristic, where the outer layer explores different resource allocations in a core-by-core manner, and the inner layer determines the task allocation under the resource allocation specified by the outer layer. To deal with the complexity of the large design space, we (i)~propose Pareto-optimality-based pruning techniques in the outer layer to reduce the exploration space by systematically eliminating suboptimal resource allocations and (ii)~model the inner-layer task allocation as a knapsack problem that can be efficiently solved by dynamic programming (DP). 

\subsection{Outer Layer: Resource Allocation Based on Pareto-Pruned Search}

\input{algorithms/outer}

Algorithm~\ref{alg:outer} describes the procedures in the outer layer. The algorithm takes the system specifications as input, including the set of real-time tasks $\tau$, the number of cores $M$, the number of memory bandwidth partitions $B$, the number of cache partitions $K$, and an algorithmic parameter $\gamma$ used for controlling the precision of the dynamic programming in the inner layer. The output of the algorithm is a set of Pareto-optimal solutions within the explored search space, where no solution dominates another in terms of all objective values. 

The algorithm employs a breadth-first search strategy to explore resource and task allocations on a core-by-core basis, with the maximum search depth equal to the total number of cores $M$. 
The search begins by initializing an empty set of complete solutions $\Omega^*$ and a single initial partial solution $\omega_0$ (lines 1-2). This initial solution has empty task and resource allocations, denoted by $\omega_0.\overrightarrow{\tau}$, $\omega_0.\overrightarrow{b}$, and $\omega_0.\overrightarrow{k}$, and retains all system resources, represented by $\omega_0.\overline{B}$ and $\omega_0.\overline{K}$ (lines 1-2). 
The set of current partial solutions $\Omega$ is initialized with the initial solution $\omega_0$ (line 3). 

During each search iteration (lines 4-18), indexed by $m$, each partial solution in $\Omega$ is extended into multiple new partial solutions, which include the task and resource allocations for the corresponding core $m$. 
Resource allocation is generated by enumerating all possible combinations of memory bandwidth and cache partitions, ranging from one partition to the number of remaining available partitions in the system (lines 7-8). 
The resource allocation of each extended partial solution is checked for dominance against the complete solutions in $\Omega^*$ using Pruning Rule~\ref{rule:complete}: 
\begin{prune}[Dominance against complete solutions]
\label{rule:complete}
    A solution $\omega_1$ is not dominated by a complete solution $\omega_2$ if and only if
    \begin{equation}
        (\omega_1.\overline{B} > \omega_2.\overline{B}) ~\lor~ (\omega_1.\overline{K} > \omega_2.\overline{K}).
    \end{equation}
\end{prune}
\noindent
If the resource usage of an extended partial solution is dominated by any complete solution in $\Omega^*$, it is pruned from further exploration, skipping the generation of its task allocation (line 9). 
For non-dominated solutions, the inner layer (Algorithm~\ref{alg:inner}) is invoked to determine the task allocation under the specified resource allocation, resulting in an extended partial solution $\omega'$ (line~10). 

Each $\omega'$ is then checked for completeness. If no tasks remain unassigend in $\omega'$, it is added to the complete solution set $\Omega^*$ (line 11 - 12). Then, $\Omega^*$ is updated to remove any solutions dominated by $\omega'$ (line~13). If there still exist remaining tasks to be assigned, the extended partial solution is checked for feasibility (line 15) using Pruning Rule~\ref{rule:feasibility}:  
\begin{prune}[Feasibility]
\label{rule:feasibility}
    A partial solution $\omega$ is infeasible if
    \begin{equation}
        \sum_{\tau_i \in \omega.\overline{\tau}}{U(\tau_i, \omega.\overline{B}, \omega.\overline{K}) > \overline{M}},
    \end{equation}
    where $\overline{M} = M - m$ is the number of remaining cores at search depth $m$. 
\end{prune}
\noindent
This rule gives a sufficient condition for infeasibility. Specifically, the total utilization of the remaining tasks, calculated under the assumption that all remaining resources are fully allocated to every core, represents a lower bound on the total remaining task utilization. If this lower bound exceeds the number of remaining cores, it indicates that the solution cannot satisfy the scheduling constraint and should be pruned from further exploration. 
It is important to note that this condition is sufficient but not necessary. In other words, while solutions that fail this check are guaranteed to be infeasible, solutions that pass the check may still turn out to be infeasible upon further exploration. 

If a solution fails the feasibility test, it will be pruned from further exploration. Otherwise, the solution is added to the new partial solution set $\Omega'$ (line 16). To maintain efficiency, we further remove dominated partial solutions in $\Omega'$ using Pruning Rule~\ref{rule:partial}:
\begin{prune}[Dominance between partial solutions]
\label{rule:partial}
    A partial solution $\omega_1$ is not dominated by another partial solution $\omega_2$ if and only if
    \begin{equation}
        (\omega_1.\overline{B} > \omega_2.\overline{B}) ~\lor~ (\omega_1.\overline{K} > \omega_2.\overline{K}) ~\lor~ (\omega_1.\overline{U} < \omega_2.\overline{U}),
    \end{equation}
    where $\omega_1.\overline{U}$ represents the remaining scheduling demand, defined as the sum of reference utilizations of the remaining tasks, \ie, $\omega_1.\overline{U} = \sum_{\tau_i \in \omega_1.\overline{\tau}}{\widehat{U}(\tau_i)}$. 
\end{prune}
\noindent
This rule intuitively prioritizes solutions with a lower remaining scheduling demand and higher remaining resources, which are more likely to be schedulable. We note that the reference utilization is not an exact indicator of scheduling demand, and thus, a dominated partial solution by this rule is not guaranteed to be worse than the non-dominated ones. However, this pruning rule effectively limits the number of partial solutions generated in each search iteration (at most $B \cdot K$, since in the worst case, for each combination of remaining memory bandwidth and cache partitions, only the partial solution with the lowest scheduling demand is kept). Evaluation results in Section~\ref{sec:eval_results} demonstrate the effectiveness of this strategy. 

The algorithm terminates when all possible resource allocations have been explored. The final set of complete solutions $\Omega^*$ contains all Pareto-optimal (non-dominated) solutions found by the algorithm.

\subsection{Inner Layer: Task Allocation Based on Dynamic Programming}

The inner layer determines the task assignment to a core with the resource allocation specified by the outer layer. 
The key idea is to assign tasks to the core to maximize the scheduling demand of assigned tasks $\tau'$ while satisfying the schedulability constraint, given the allocated resources ($b$ memory bandwidth and $k$ cache partitions). 
The scheduling demand of a task $\tau_i$ is measured by its reference utilization $\widehat{U}(\tau_i)$, and the schedulability constraint is expressed as $\sum_{\tau_i \in \tau'}{U(\tau_i, b, k)} \leq 1$. 
This problem can be modeled as a 0-1 knapsack problem~\cite{martello1990knapsack}, where:
\begin{itemize}
    \item The remaining tasks to be assigend $\omega.\overline{\tau}$ are regarded as items to be packed into the knapsack;
    \item The value of each item is the task's reference utilization $\widehat{U}(\tau_i)$;
    \item The size of each item is the task's utilization under the current resource allocation $U(\tau_i, b, k)$;
    \item The knapsack's capacity corresponds to the utilization bound of the EDF scheduling policy, which is 1. 
\end{itemize}
The knapsack problem is solved by a dynamic programming (DP) approach~\cite{martello1990knapsack}, as outlined in Algorithm~\ref{alg:inner}. 
\input{algorithms/inner}

The algorithm takes the partial solution to be extended $\omega$, along with the allocated memory bandwidth $b$ and cache partitions $k$, and a scaling parameter $\gamma$. Since dynamic programming requires integer item sizes, the scaling parameter $\gamma$ is introduced to inflate fractional task utilizations and the utilization bound to integer values. In particular, each task's utilization is scaled and rounded up using the ceiling function to ensure a safe utilization bound (line 5). 
At the beginning of the algorithm, we initialize the new partial solution $\omega'$ by appending allocated resources ($b$ and $k$) to the corresponding resource allocation vectors and recalculating the remaining resources available in the system (lines 1-2). 
Then, a 2D DP array $dp[0...n,0...\gamma]$ is initialized, where $n$ is the number of unassigned tasks. The value in $dp[i,j]$ stores the maximum total reference utilization achievable for the first $i$ tasks with a scaled utilization bound $j$. 
Then, the DP array is updated iteratively. If the task's scaled utilization exceeds the current bound $j$, the task is excluded (line 8); otherwise, the maximum reference utilization is determined by taking the better of two options: excluding the task or including it and adding its reference utilization (line 9). 
After constructing the DP array, the selected tasks are retrieved by backtracking through the array (lines 13–17).
Finally, the algorithm appends the assigned set of tasks to the solution's task allocation vector and returns the extended partial solution $\omega'$ (lines 17-18).

\subsection{Time Complexity}
\label{sec:complexity}

In the outer layer, the maximum number of search iterations is given by $M$, and the number of partial solutions to be extended in each iteration is bounded by $B \cdot K$ under the dominance relationship defined in Pruning Rule~\ref{rule:partial}. In the worst case, all possible resource allocations are considered for each partial solution, resulting in $\mathcal{O}(B^2 K^2)$ invocations of the inner layer to determine task allocations. 
Since the DP approach has a time complexity of $\mathcal{O}(N \gamma)$, the overall complexity of the algorithm is $\mathcal{O}(N B^2 K^2 M \gamma)$. 
In this paper, we did preliminary experiments with different $\gamma$ values ranging from 1000 to 10000. As the results are consistent across these values, we only present the results with $\gamma=1000$.

%% file: algorithms/outer.tex
{
\SetInd{0.2em}{2em}
\begin{algorithm} [t]
    \caption{Multi-Objective Search with Pareto-Based Pruning (Outer Layer)}
    \label{alg:outer}
    
    \KwIn{$\tau$, $M$, $B$, $K$, $\gamma$\;}
    \KwOut{$\Omega^*$: The non-dominated solutions\;}
    
    $\Omega^* \gets \varnothing$\;
    
    Initialize $\omega_0$ with empty allocation vectors: $\omega_0.\overrightarrow{\tau}, \omega_0.\overrightarrow{b}, \omega_0.\overrightarrow{k}$, and full remaining tasks and resources: $\omega_0.\overline{\tau} \gets \tau, \omega_0.\overline{B} \gets B, \omega_0.\overline{K} \gets K$\;

    
    $\Omega \gets \{\omega_0\}$\;
    
	\For{$m \gets 1$ \KwTo $M$}
    {
        $\Omega' \gets \varnothing$\;
	    \For{$\omega \in \Omega$}
        {
            \For{$b \gets 1$ \KwTo $\omega.\overline{B}$}
            {
                \For{$k \gets 1$ \KwTo $\omega.\overline{K}$}
                {
                    \If(\tcp*[f]{Rule 1}){$(\omega.\overline{B}{-}b, \omega.\overline{K}{-}k)$ not dominated by $\Omega^*$}
                    {
                        $\omega' \gets \textit{TaskAlloc}(\omega.\overline{\tau}, b, k, \gamma)$\;
                        \eIf{$\omega'.\overline{\tau} = \varnothing$}
                        {
                            $\Omega^* \gets \Omega^* \cup \{\omega'\}$\;
                            \textit{PruneCompleteSols}($\Omega^*$); \tcp*[f]{Rule 1}
                        }
                        {
                            \If(\tcp*[f]{Rule 2}){\textit{IsFeasible}$(\omega')$}
                            {
                                $\Omega' \gets \Omega' \cup \{\omega'\}$\;
                                \textit{PrunePartialSols}($\Omega'$); \tcp*[f]{Rule 3}
                            }
                        }
                    }
                }
            }
	    }
     $\Omega \gets \Omega'$\;
	}
	\KwRet{$\Omega^*$}\;
	
\end{algorithm}
}

%% file: algorithms/inner.tex
{
\SetInd{0.2em}{2em}
\begin{algorithm} [t]
    \caption{Dynamic Programming Based Task Allocation (Inner Layer)}
    \label{alg:inner}
    
    \KwIn{$\omega$, $b$, $k$, $\gamma$\;}
    \KwOut{$\omega'$: The extended solution\;}

    $\omega'.\overrightarrow{\tau} \gets \omega.\overrightarrow{\tau}$; Append $b$ to $\omega'.\overrightarrow{b}$; Append $k$ to $\omega'.\overrightarrow{k}$; $n \gets |\omega.\overline{\tau}|$\;

    $\omega'.\overline{b} \gets \omega.\overline{b} - b$; $\omega'.\overline{k} \gets \omega.\overline{k} - k$\;

    Initialize a 2D array $dp[0 \dots n,0 \dots \gamma]$ with zeros\;

    \For{$i = 1$ \KwTo $n$}
    {
        $u \gets \lceil U(\omega.\overline{\tau}_{i-1}, b, k) \cdot \gamma \rceil$\;
        $\widehat{u} \gets \widehat{U}(\omega.\overline{\tau}_{i-1})$\;
        \For{$j = 0$ \KwTo $\gamma$}
        {
            \eIf{$u > j$}
            {
                $dp[i,j] \gets dp[i - 1,j]$;
            }
            {
                $dp[i,j] \gets \max(dp[i - 1,j], dp[i - 1,j - u] + \widehat{u})$\;
            }
        }
    }

    $\tau' \gets \varnothing$;
    $\textit{j} \gets \gamma$\;
    
    \For{$i = n$ \KwTo $1$}
    {
        \If{$dp[i,j] \neq dp[i - 1,j]$}
        {
            $\tau' \gets \tau' \cup \{\omega.\overline{\tau}_{i-1}\}$\;
            $\textit{j} \gets \textit{j} - \lceil U(\omega.\overline{\tau}_{i-1}, b, k) \cdot \gamma \rceil$\;
        }
    }
    Append $\tau'$ to $\omega'.\overrightarrow{\tau}$\;

	\KwRet{$\omega'$}\; 
	
\end{algorithm}
}

%% file: sec_tools.tex
\section{Real-World Benchmarking with Resource Partitioning}
\label{sec:tools}

\subsection{Resource Partitioning Configuration}

We profile real-world benchmarks on the embedded AMD UltraScale+ ZCU102~\cite{zcu-102} platform, equipped with four Cortex-A53 cores. Each core has a private 32~KB, 4-way set-associative L1 data cache\footnote{Each core also has a private 2-way set-associative L1 instruction cache, but this paper focuses on the data cache.} and shares a 1~MB, 16-way set-associative last-level cache (LLC).
We use benchmarks from two test suites: 43 benchmarks from \sdvbshort~\cite{venkata2009sd, rt-bench} and 20 benchmarks from \parsec~\cite{parsec}.\footnote{We do not use all the benchmarks from the test suites because the remaining benchmarks either fail to be built or take too long to execute on our platform.} Both test suites have also been previously adopted in similar problem settings~\cite{XPCLLLL:19, sun2023co}. 
Notably, \parsec was ported to the ARM 64-bit architecture to ensure compatibility with our platform. 

The profiles were generated by systematically limiting the availability of the LLC and memory bandwidth resources allocated to the benchmarks.
This was achieved transparently through the Jailhouse hypervisor~\cite{jailhouse, minerva-jailhouse}, which provides mechanisms for cache set partitioning and memory bandwidth control.\footnote{We use the publicly available Minervasys (\url{https://github.com/Minervasys/Jailhouse}) fork of Jailhouse that implements cache and bandwidth partitioning.}
Specifically, the hypervisor
employs coloring techniques~\cite{KSMCV:19} for cache set partitioning and integrates the MemGuard~\cite{yun2013memguard, STPMBZC:21} to regulate bandwidth allocation. 
The benchmarks were executed on a single core of the ZCU102 under a virtualized Linux operating system (Kernel 6.1.30).
Slowdown profiles were computed relative to configurations with full hardware resource availability (\ie, the Jailhouse configuration gives full hardware resources to the virtual machine running Linux). 
Each configuration was repeated 10 times, with a total runtime cap of 10 minutes per benchmark. 

The granularity of cache set partitioning depends on the architecture of ZCU102.
With a 1~MB, 16-way set associative L2 and 64-byte cache lines (6~offset bits), the cache has 1024 sets (10~index bits).
Given a 4 KB page size, which aligns to the last 12 address bits, four bits in the index range are available for partitioning. This allows for a maximum of 16 partitions, each corresponding to 64 KB of the L2 cache. 
While this scheme offers high granularity, it can lead to underutilization of the L1 data cache due to its associativity, as up to half of the L1 cache might remain unused in certain configurations. 
To evaluate this trade-off, we chose to enable all possible coloring configurations, prioritizing higher granularity over limiting the maximum number of partitions to avoid L1 underutilization. We experimentally validated the impact of this choice by running a dedicated set of experiments on \sdvbshort benchmarks with various cache-partitioning configurations, including those that colored the L1 cache. For each configuration, we measured execution time while keeping the number of cache partitions constant.
The results showed that across all benchmarks, the variation in execution time between configurations was under 5\%. For larger data sets, such as \texttt{VGA}, this variation dropped below 1\%, indicating that the penalty of L1 underutilization is minimal in practice.
Therefore, we use 16 cache partitions each of 64 KB size of L2 in our benchmarking.

Bandwidth partitioning was similarly controlled using MemGuard in the Jailhouse hypervisor.
We validated a MemGuard period of 1 ms by evaluating its overheads and benchmark execution times.
This value aligns with the findings from previous studies \eg~\cite{yun2013memguard, MemPol}.
To establish the maximum usable budget, we increased the allocated bandwidth in steps of 64 MB/s until the execution time for each benchmark stabilized, with variations below 5\%. 
This approach yielded a total of 15 distinct equally-sized bandwidth partitions, with a maximum guaranteed bandwidth of 960 MB/s.
We note that the approach and the obtained bandwidth values are consistent with prior works, \eg, \cite{yun2013memguard, MemPol, XPCLLLL:19}.

\subsection{Profiling Results}

Figure~\ref{fig:bench_both} reports the obtained slowdown profiles (relative to the execution without partitioning) for representative benchmarks from the \sdvbshort and \parsec test suites.%
\footnote{The complete profiling results are available in the supplemental material of the submission and will be made publicly available as a technical report upon the acceptance of the paper.}
In the figures, each column corresponds to the slowdowns of a single benchmark, while different rows within the same column show results for varying input sizes. The color scale (logarithmic) indicates the slowdown compared to execution without any partitioning: lighter shades denote more significant slowdowns, whereas darker shades represent near-ideal performance. 
The benchmarks have been chosen to illustrate representative slowdown trends, but we note that our evaluation presented in Section~\ref{sec:exp} is based on the full set of benchmarks, demonstrating how our co-allocation strategies handle both extreme and intermediate resource sensitivities. 

\begin{figure}
    \centering
    \input{figures/bench_sdvbs}
    \input{figures/bench_parsec}
    \caption{Slowdown profiles for representative benchmarks from different test suites. Columns present the results for one benchmark for different inputs (rows). Numerical slowdown values for the extreme configurations are indicated in the corners of each plot.}
    \label{fig:bench_both}
\end{figure}

From the \sdvbshort suite (Figure~\ref{fig:big_bench_sdvbs}), we observe that the \texttt{disparity} and \texttt{mser} benchmarks are highly sensitive to changes in both memory and cache allocations. When these resources are scarce, slowdowns exceeding 20$\times$ can occur, but they decrease steadily as additional resources are assigned.  
Notably, \texttt{disparity-vga} and \texttt{mser-vga} show a significant slowdown of 11.3$\times$ and 8.7$\times$ respectively when only minimal memory bandwidth is available.
The \texttt{stich} benchmark exhibits less uniform behavior. For most inputs, its performance remains stable across a wide range of cache and memory bandwidth settings. However, the \texttt{stich-vga} input is a notable exception: increasing the number of cache partitions sharply reduces the slowdown from 34.3$\times$ to 1.2$\times$. 
\texttt{sift} is generally less sensitive to memory bandwidth and mostly insensitive to the number of cache partitions (except for $k=1$ with inputs \texttt{vga} and \texttt{cif}, where the workload cannot finish execution within the allotted time). 
In fact, with only one memory bandwidth partition assigned, there is a minor slowdown difference between one cache partition and a fully allocated cache, and it is sufficient to assign a few additional memory bandwidth partitions to bring the slowdown to 1.0.

Compared to \sdvbshort, the benchmarks from the \parsec suite (in Figure~\ref{fig:big_bench_parsec}) show comparatively milder slowdowns overall, typically within the 1–13$\times$ range. 
Notably, in all \parsec benchmarks, one single cache partition is not sufficient for the workload to be completed within the allotted time. This is reflected in Figure~\ref{fig:big_bench_parsec}, where results start from a minimum of two cache partitions.
Among these benchmarks, \texttt{canneal} reaches its highest slowdown of around 13$\times$ when both cache and memory bandwidth are minimized. \texttt{bodytrack} and \texttt{ferret} are moderately sensitive to resource allocations, showing gradual improvement as more resources become available. \texttt{freqmine} is among the least sensitive in this set, seldom exceeding a 2–3$\times$ slowdown even under constrained conditions.

%% file: figures/bench_sdvbs.tex
\begin{subfigure}{\linewidth}
    \centering
    \begin{tikzpicture}[scale=0.96]
    \small
        \node[] at (0,0) {\includegraphics[width=3cm]{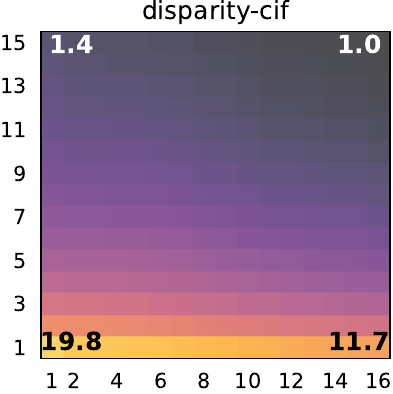}};
        \node[] at (0,3.2) {\includegraphics[width=3cm]{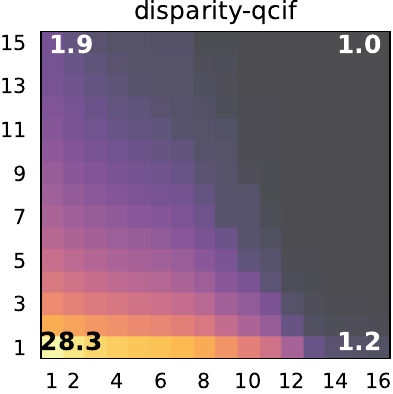}};
        \node[] at (0,6.4) {\includegraphics[width=3cm]{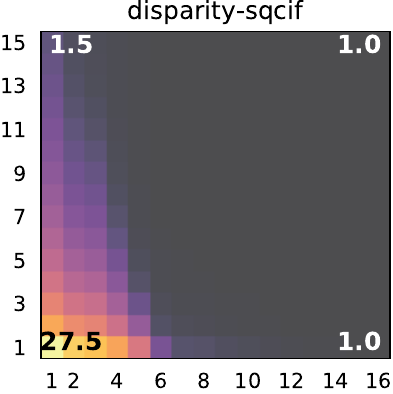}};
        \node[] at (0,9.6) {\includegraphics[width=3cm]{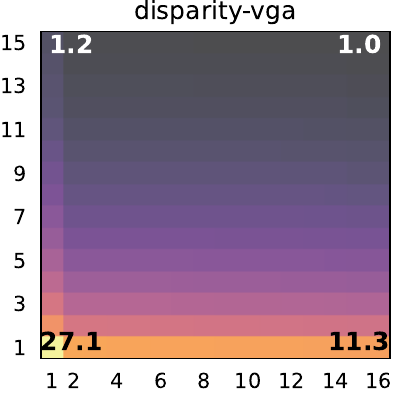}};
        
        \node[] at (3,0) {\includegraphics[width=3cm]{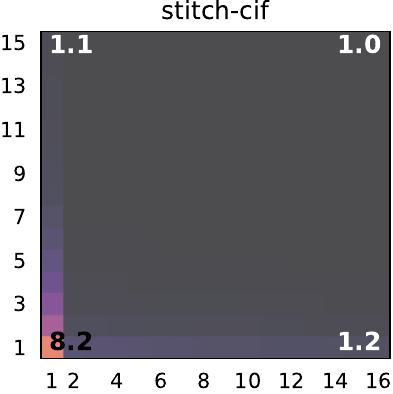}};
        \node[] at (3,3.2) {\includegraphics[width=3cm]{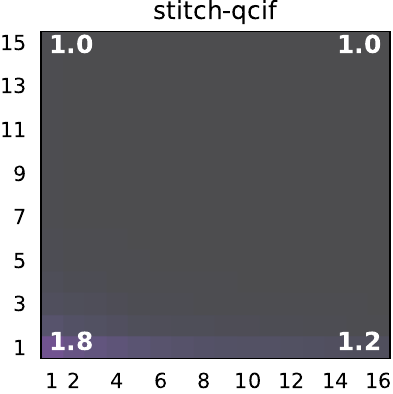}};
        \node[] at (3,6.4) {\includegraphics[width=3cm]{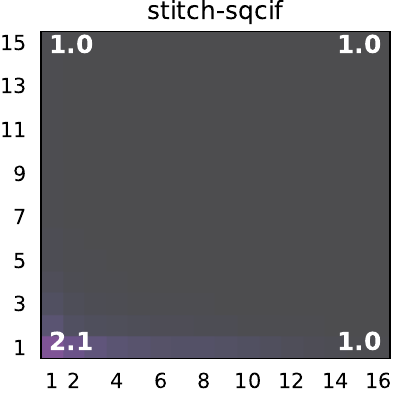}};
        \node[] at (3,9.6) {\includegraphics[width=3cm]{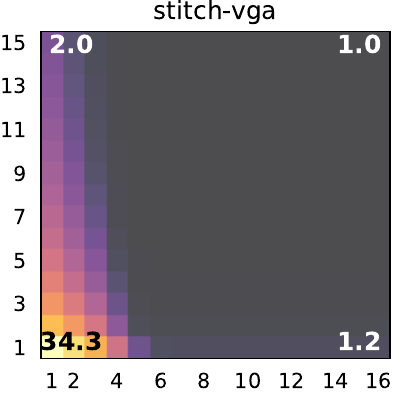}};
        
        \node[] at (6,0) {\includegraphics[width=3cm]{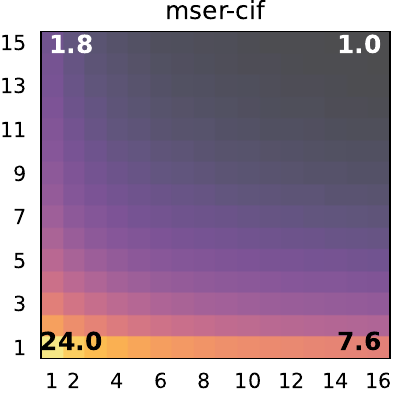}};
        \node[] at (6,3.2) {\includegraphics[width=3cm]{figures/sdvbs/mser-qcif.pdf}};
        \node[] at (6,6.4) {\includegraphics[width=3cm]{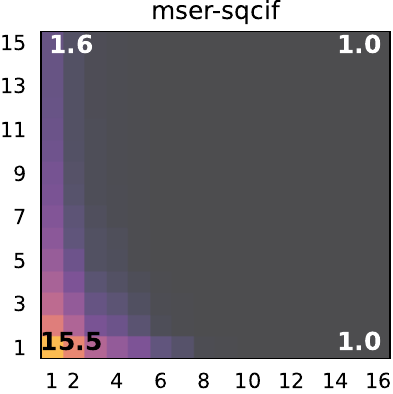}};
        \node[] at (6,9.6) {\includegraphics[width=3cm]{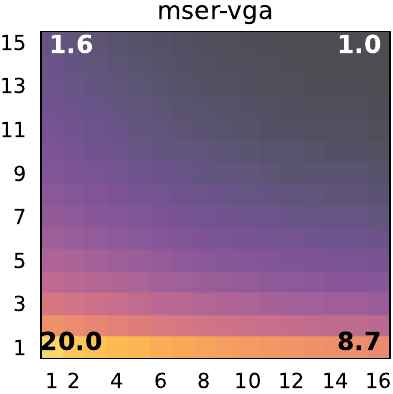}};
        
        \node[] at (9,0) {\includegraphics[width=3cm]{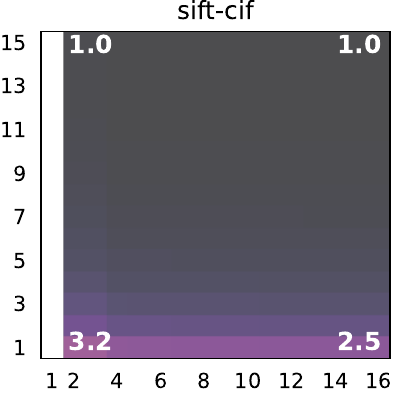}};
        \node[] at (9,3.2) {\includegraphics[width=3cm]{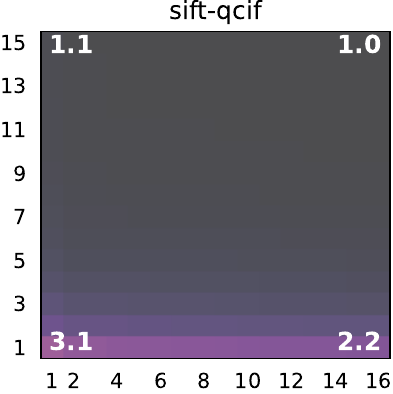}};
        \node[] at (9,6.4) {\includegraphics[width=3cm]{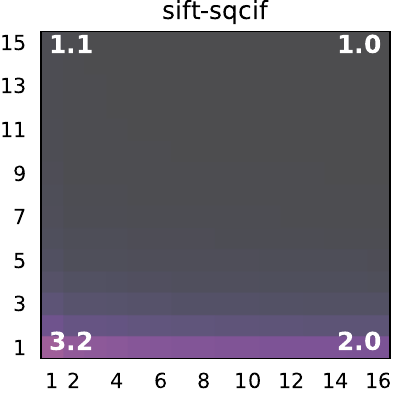}};
        \node[] at (9,9.6) {\includegraphics[width=3cm]{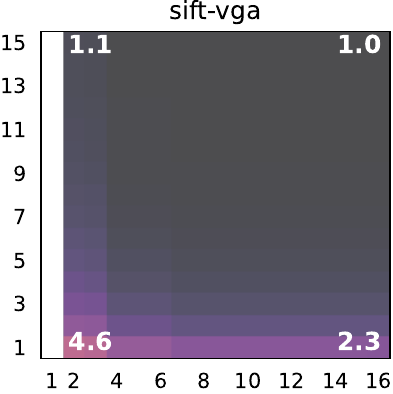}};       
        \node[] at (11,5) {\includegraphics[height=5.8cm]{figures/sdvbs/disparity-sqcif_colorbar.pdf}};        
        

        \node[rotate=90] at (-1.8, 4.5) {Number of Memory Bandwidth Partitions};
        \node[] at (4.5, -1.8) {Number of Cache Partitions};
        \node[rotate=90] at (11.5, 5) {Slowdown (log-scale)};
    \end{tikzpicture}
    \caption{Representative slowdown profiles for the \sdvbshort test suite.}
    \label{fig:big_bench_sdvbs}
\end{subfigure}

%% file: figures/bench_parsec.tex
\begin{subfigure}{\linewidth}
    \centering
    \begin{tikzpicture}[scale=0.96]
    \tiny
        \node[] at (0,0) {\includegraphics[width=3cm]{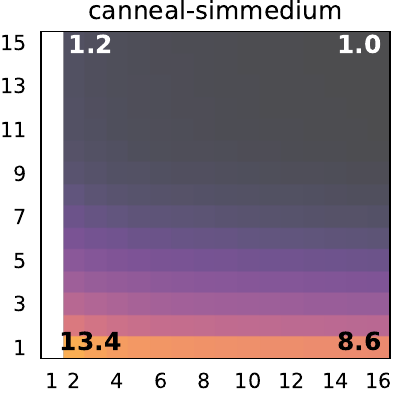}};
        \node[] at (0,3.2) {\includegraphics[width=3cm]{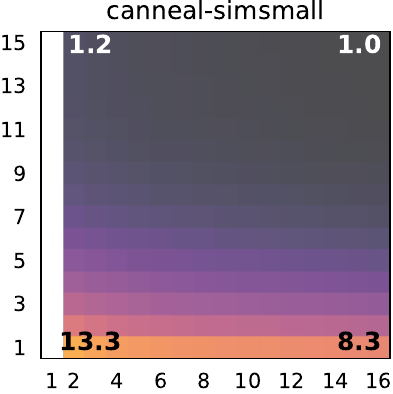}};
        
        \node[] at (3,0) {\includegraphics[width=3cm]{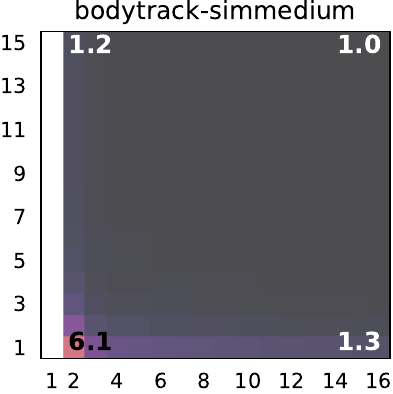}};
        \node[] at (3,3.2) {\includegraphics[width=3cm]{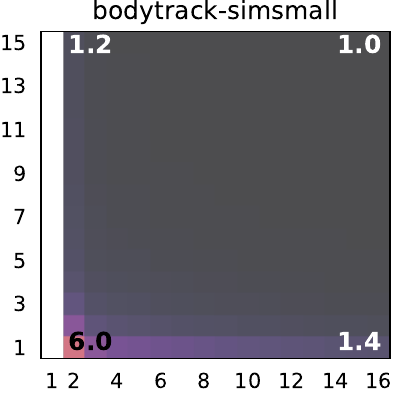}};
        
        \node[] at (6,0) {\includegraphics[width=3cm]{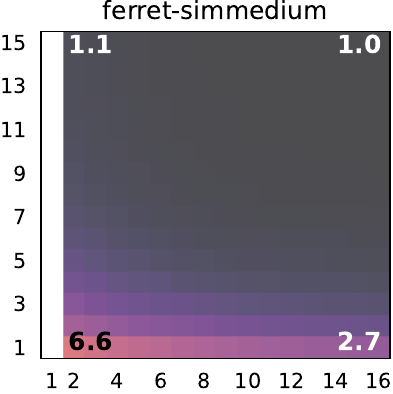}};
        \node[] at (6,3.2) {\includegraphics[width=3cm]{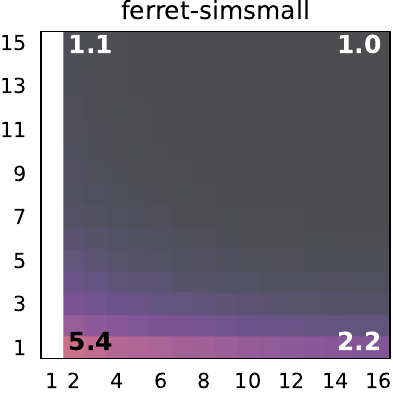}};
        
        \node[] at (9,0) {\includegraphics[width=3cm]{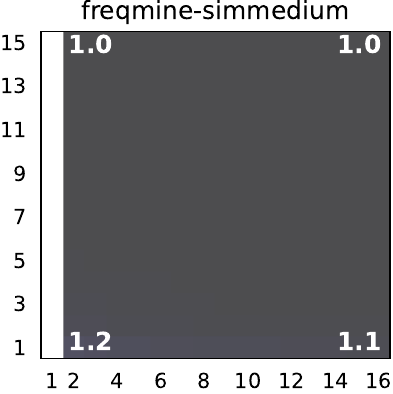}};
        \node[] at (9,3.2) {\includegraphics[width=3cm]{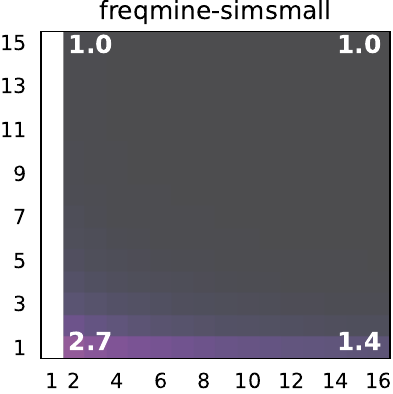}};

        \node[] at (11,1.5) {\includegraphics[height=5.8cm]{figures/sdvbs/disparity-sqcif_colorbar.pdf}};    

        \node[rotate=90] at (-1.8, 1.5) {Number of Memory Bandwidth Partitions};
        \node[] at (4.5, -1.8) {Number of Cache Partitions};
        \node[rotate=90] at (11.5, 1.5) {Slowdown (log-scale)};
    \end{tikzpicture}
    \caption{Representative slowdown profiles for the \parsec test suite.}
    \label{fig:big_bench_parsec}
\end{subfigure}

%% file: sec_experiments.tex
\section{Performance Evaluation}
\label{sec:exp}

\subsection{Evaluation Setup}
\label{sec:eval_setup}

\subsubsection{Task Sets}
We evaluate the proposed algorithm using real-time task sets generated based on the benchmarks profiled in Section~\ref{sec:tools}. 
The total number of cores and resources available in the system is set according to the real-world embedded hardware used in our profiling, \ie, $M=4$, $B=15$, and $K=16$. 
We generate task sets with size $N \in \{20, 40, 60\}$ and consider different task set reference utilizations $\widehat{U} = \sum_{\tau_i \in \tau}{\widehat{U}(\tau_i)}$ ranging from $1.0$ to $M$, with a step of $0.1$. For each combination of $N$ and $\widehat{U}$, we generate 100 random task sets per benchmark suite, resulting in a total of 18,600 task sets. 
For each task set, we use the \texttt{DRS} generator~\cite{griffin2020generating} to assign task reference utilizations, ensuring their sum matches the target $\widehat{U}$. For each task in a task set, a slowdown profile is randomly sampled from the benchmark programs and used to calculate the task utilizations under different resource allocations.

\subsubsection{Baseline Algorithms and Implementation}
We compare the proposed MMO framework with the following baselines:
\begin{itemize}
    \item CaM~\cite{XPCLLLL:19}: the state-of-the-art heuristic that addresses the same task-resource co-allocation problem;
    \item ILP-B, ILP-K: the single-objective 0-1 ILPs to minimize memory bandwidth and cache usage, respectively, presented in Section~\ref{sec:sys_model}; 
    \item MIP-B, MIP-K: the single-objective mixed-integer programs (MIPs) to minimize memory bandwidth and cache usage, respectively, developed in the extended version\footnote{\url{https://www.cis.upenn.edu/~linhphan/papers/rtas19-CaM-techreport.pdf}} of \cite{XPCLLLL:19}. 
\end{itemize} 
The 0-1 ILPs and MIPs are solved by a standard mathematical programming solver, Gurobi~\cite{gurobi}. 
Given the scalability challenges of mathematical programming solvers for large decision spaces, we set a one-hour time limit per task set. If the solver fails to prove optimality within this limit, the best solution found is used for evaluation.

All algorithms are implemented in Python 3.12, with Gurobi version 11.0. Experiments were conducted on a workstation equipped with AMD EPYC 7763 CPUs running GNU/Linux.

\subsection{Evaluation Results}
\label{sec:eval_results}

We evaluate the algorithms on the 18,600 task sets and compare them in terms of schedulability ratio, memory bandwidth usage, cache usage, and solving efficiency. Additionally, we report the number of Pareto-optimal solutions obtained by MMO.

\subsubsection{Schedulability}
Figure~\ref{fig:sched} presents schedulability ratio \textit{w.r.t.} reference utilization across all tested task set sizes ($N\in\{20,40,60\}$).
MMO consistently achieves the highest schedulability ratio across all utilizations and benchmark suites. In particular, MMO identifies up to 62.67\% and 50.34\% more schedulable task sets than CaM on the \textit{PARSEC} and \textit{SD-VBS} benchmark suites, respectively. The results validate the effectiveness of MMO's design improvement over CaM, as discussed in Section~\ref{sec:related_work_multi_resource_co_alloc}. 
Among the mathematical programming models, ILP-B and ILP-K outperform MIP-B and MIP-K because the linear constraint formulation in the developed 0-1 ILP is solved more efficiently by the solver than the non-linear constraints in the MIP. 
However, they still fail to obtain the optimal solutions within one hour for many task sets, where the best-found solutions are used for evaluation.

\input{figures/sched}

\input{figures/bw}

\input{figures/cp}

\input{figures/sols}

\subsubsection{Resource Usage}
Figure~\ref{fig:bw} and Figure~\ref{fig:cp} illustrate the memory bandwidth and cache usage, respectively, measured in the number of partitions required for scheduling. If a task set is found to be unschedulable by an algorithm, the resource usage is considered as the total available partitions of that resource in the system. 
The y-axis in each figure represents the average resource usage across all evaluated task sets at a given system utilization level (x-axis). At each utilization level, the solutions with the lowest memory bandwidth and cache usage among the Pareto-optimal set are selected for comparison in Figures~\ref{fig:bw} and~\ref{fig:cp}, respectively. 
Results show that MMO effectively schedules the task sets with the fewest average resources compared to the baselines. 
For memory bandwidth and cache usage, MMO achieves similar performance to ILP-B and ILP-K, the single objective derivatives of the 0-1 ILP, respectively. 
CaM outperforms MIP-K/ILP-K and MIP-B/ILP-B for memory bandwidth and cache usage, respectively, because these mathematical programs are not configured to minimize the corresponding resource usage in their objective functions. However, all these models achieve better performance than CaM for the objectives they optimize.

\subsubsection{Number of Non-Dominated Solutions}
Unlike the baselines, MMO can generate multiple solutions with non-dominated objective values. Figure~\ref{fig:sols_hist} shows the histogram of the number of non-dominated solutions found by MMO for all schedulable task sets. Results show that MMO finds at least 3 non-dominated solutions for over 90\% of schedulable task sets. 
Figure~\ref{fig:sols_example} provides an example Pareto front obtained by MMO and the best solutions found by the baseline algorithms for a task set with $N=40$ and $\widehat{U} = 2.4$. 
The Pareto front is comprised of 7 non-dominated solutions, with memory bandwidth usage ranging from 6 to 15 partitions and cache usage ranging from 6 to 14 partitions. 
For the baselines, ILP-B and MIP-B obtain the same solution with memory bandwidth usage equal to the minimum one found by MMO, but their cache usage is larger and therefore dominated by MMO's Pareto front. 
In terms of cache usage, ILP-K yields a solution that requires one additional cache partition compared to MMO’s minimum, both with the same memory bandwidth requirement. 
Although ILP-K is designed to minimize cache usage, the solution obtained within the one-hour time limit is suboptimal, as the solver did not converge to the optimal solution within that time.
The non-dominated solutions in MMO's Pareto front offer system designers the flexibility to select the most suitable solution based on their preferences.

\subsubsection{Solving Efficiency}
Table~\ref{tab:efficiency} summarizes the runtime performance of each algorithm concerning different task set sizes. MMO achieves the lowest average running time among all algorithms in comparison. Notably, the average solving time of MMO is around 3x lower than CaM. The maximum running time of MMO across all task sets is 17.1 seconds, significantly lower than the baselines.

\input{tables/efficiency}

%% file: figures/sched.tex
\begin{figure}[t]
    \centering
    \begin{subfigure}{\linewidth}
      \centering
      \includegraphics[width=0.9\linewidth]{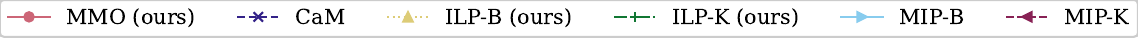}
    \end{subfigure}
    \begin{subfigure}{.49\textwidth}
      \centering
      \includegraphics[width=0.95\linewidth]{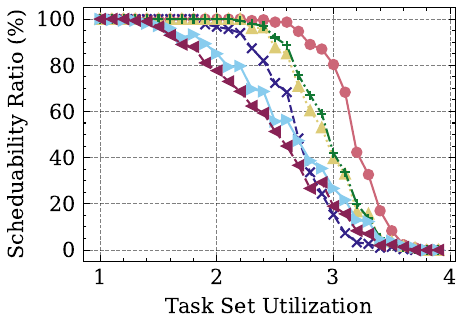}
      \caption{PARSEC Benchmarks.}
      \label{fig:sched_parsec}
    \end{subfigure}%
    \begin{subfigure}{.49\textwidth}
      \centering
      \includegraphics[width=0.95\linewidth]{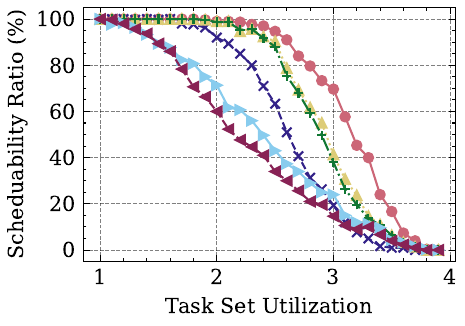}
      \caption{SD-VBS Benchmarks.}
      \label{fig:sched_sdvbs}
    \end{subfigure}%
    \caption{Schedulability ratio (\textit{i.e.}, $\frac{\text{\#schedulable task sets}}{\text{\#total task sets}} \times 100$ \%).}
    \label{fig:sched}
\end{figure}

%% file: figures/bw.tex
\begin{figure}[t]
    \centering
    \begin{subfigure}{\linewidth}
      \centering
      \includegraphics[width=0.9\linewidth]{figures/parsec/sched_u_legend.pdf}
    \end{subfigure}
    \begin{subfigure}{.49\textwidth}
      \centering
      \includegraphics[width=0.95\linewidth]{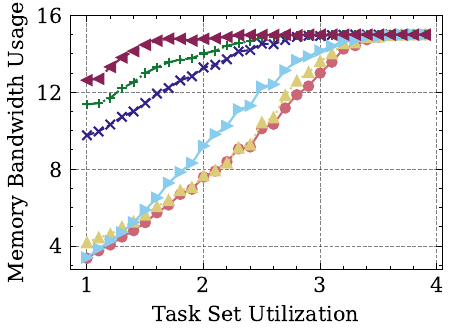}
      \caption{PARSEC Benchmarks.}
      \label{fig:bw_parsec}
    \end{subfigure}%
    \begin{subfigure}{.49\textwidth}
      \centering
      \includegraphics[width=0.95\linewidth]{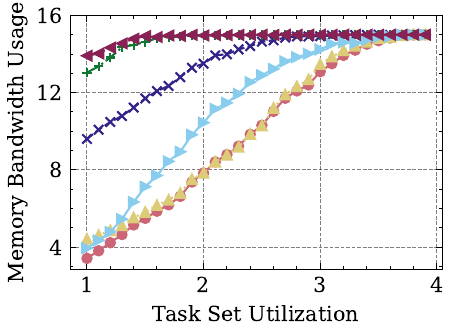}
      \caption{SD-VBS Benchmarks.}
      \label{fig:bw_sdvbs}
    \end{subfigure}%
    \caption{Memory bandwidth usage (number of partitions).}
    \label{fig:bw}
\end{figure}

%% file: figures/cp.tex
\begin{figure}[t]
    \centering
    \begin{subfigure}{\linewidth}
      \centering
      \includegraphics[width=0.9\linewidth]{figures/parsec/sched_u_legend.pdf}
    \end{subfigure}
    \begin{subfigure}{.49\textwidth}
      \centering
      \includegraphics[width=0.95\linewidth]{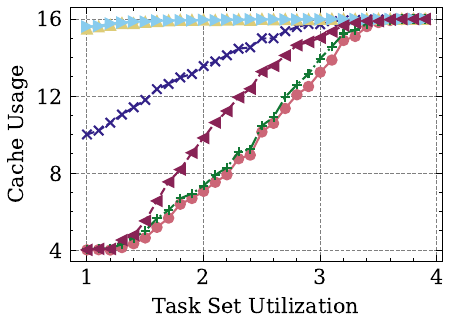}
      \caption{PARSEC Benchmarks.}
      \label{fig:cp_parsec}
    \end{subfigure}%
    \begin{subfigure}{.49\textwidth}
      \centering
      \includegraphics[width=0.95\linewidth]{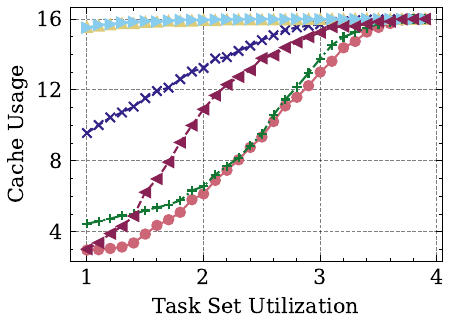}
      \caption{SD-VBS Benchmarks.}
      \label{fig:cp_sdvbs}
    \end{subfigure}%
    \caption{Cache usage (number of partitions).}
    \label{fig:cp}
\end{figure}

%% file: figures/sols.tex
\begin{figure}[t]
    \centering
    \begin{subfigure}{.49\textwidth}
      \centering
      \includegraphics[width=0.95\linewidth]{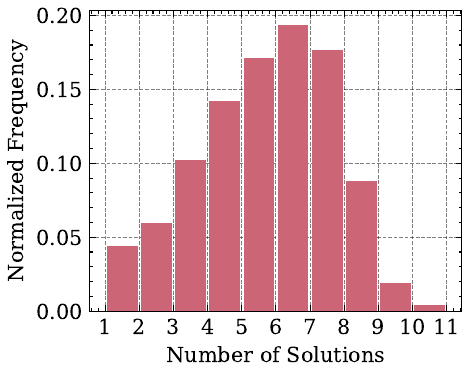}
      \caption{Histogram of MMO's non-dominated solutions.}
      \label{fig:sols_hist}
    \end{subfigure}%
    \begin{subfigure}{.49\textwidth}
      \centering
      \includegraphics[width=0.95\linewidth]{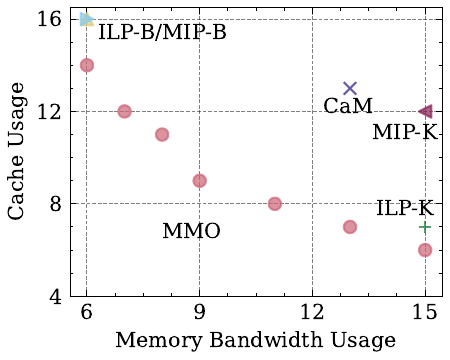}
      \caption{An example of MMO's Pareto front.}
      \label{fig:sols_example}
    \end{subfigure}%
    \caption{Pareto-optimal solutions generated by the comparison algorithms.}
    \label{fig:sols}
\end{figure}

%% file: tables/efficiency.tex
\begin{table}[ht]
\centering
\setlength{\tabcolsep}{3pt}
\begin{tabular}{c|ccc|ccc|ccc}
    \hline
    \multirow{2}{*}{Algorithm} & \multicolumn{3}{c|}{$N=20$} & \multicolumn{3}{c|}{$N=40$} & \multicolumn{3}{c}{$N=60$} \\
    \cline{2-10}
              & Avg  & Min  & Max   & Avg   & Min  & Max    & Avg   & Min  & Max    \\
    \hline
    MMO   & 2.16  & 0.04 & 8.73  & 3.75  & 0.06 & 12.0  & 5.06  & 0.08 & 17.1  \\
    CaM   & 7.94  & 0.02 & 45.6 & 12.3 & 0.09 & 67.3  & 16.8 & 0.14 & 91.6  \\
    ILP-B & 1294  & 8.43 & OOT   & 2536  & 31.2 & OOT   & 2878  & 58.9 & OOT    \\
    ILP-K & 1185  & 5.67 & OOT   & 2270  & 21.7 & OOT   & 2469  & 51.3 & OOT    \\
    MIP-B & 1756  & 1.73 & OOT   & 3390  & 21.1 & OOT   & 3485  & 34.4 & OOT    \\
    MIP-K & 1881  & 1.17 & OOT   & 3019  & 5.29 & OOT   & 3205  & 10.5 & OOT    \\
    \hline
\end{tabular}
\caption{Solving time comparison. Each entry gives the average, minimum, and maximum running times in seconds of an algorithm considering all task sets with a certain size ($N$). ``OOT'' represents ``out of time,'' meaning the solver reaches the one-hour limit.} 
\label{tab:efficiency}
\end{table}

%% file: sec_conclusion.tex
\section{Conclusion}
\label{sec:conclusion}

This paper studies a multi-objective task-resource co-allocation problem to minimize cache and memory bandwidth usage while ensuring task schedulability under a partitioned preemptive EDF scheduling policy. 
To this end, we first formulate a zero-one linear programming model (0-1 ILP), which can be solved to optimize every single objective separately using a standard mathematical solver.
To overcome the problem of exponential time complexity of the ILP and its inability to generate multiple non-dominated solutions within a reasonable time, we propose MMO, a multi-objective multi-layer optimization heuristic that optimizes both objectives simultaneously. Our approach adopts a multi-layer optimization framework where the outer layer explores resource allocations using Pareto-pruned search, while the inner one optimizes task allocation. 

We evaluate our approach by deriving slowdown profiles for real-world benchmarks from the \sdvbshort and \parsec suites executed on a real embedded AMD UltraScale+ ZCU102 platform.
We configure fine-grained resource partitioning on this MPSoC using the Jailhouse hypervisor, leveraging cache set partitioning and Memguard for memory bandwidth regulation. 
Our experiments based on these real-world benchmarks demonstrate that
our heuristics perform consistently better than the state-of-the-art in terms of schedulability, resource usage, number of non-dominated solutions, and solving efficiency.

In future work, we will investigate task-resource co-allocation strategies for parallel real-time tasks with data dependencies (\eg, modeled by directed acyclic graphs).

%% file: supplemental_material.tex
\begin{appendices}







\section*{Supplemental Material}

In this supplemental material, we present the slowdown plots for all 43 benchmarks profiled in the SD-VBS suite and all 20 benchmarks profiled in the PARSEC suite. These profiles are generated by varying the number of memory bandwidth partitions and cache partitions. Slowdown is calculated as the ratio of the benchmark's execution time with the tested memory bandwidth and cache allocation to its execution time with full memory bandwidth and cache allocation. The slowdown is plotted on a log-scale, with the x-axis representing the number of memory bandwidth partitions and the y-axis representing the number of cache partitions. Slowdown values under extreme memory bandwidth and cache partition allocations are annotated at the corners of each plot. 

Some benchmarks cannot be completed with certain cache allocations (e.g., with only one cache partition), and their corresponding slowdown values are left blank (white in the heatmap). The benchmarks are ordered alphabetically.

\section{SD-VBS Benchmarks}

    \begin{tikzpicture}[] 
    \node[] at (1,1) {\includegraphics[width=5cm]{figures/sdvbs/disparity-cif.pdf}};
    \node[rotate=90, align=center] at (-2, 1) {Number of Memory \\ Bandwidth Partitions};
    \node[] at (1.2, -1.8) {Number of Cache Partitions};
    \end{tikzpicture}
    \begin{tikzpicture}[]
    \node[] at (1,1) {\includegraphics[width=5cm]{figures/sdvbs/disparity-qcif.pdf}};
    \node[] at (4.5,1) {\includegraphics[height=5cm]{figures/sdvbs/disparity-sqcif_colorbar.pdf}};        
    \node[rotate=90, align=center] at (-2, 1) {Number of Memory \\ Bandwidth Partitions};
    \node[] at (1.2, -1.8) {Number of Cache Partitions};
    \node[rotate=90] at (5, 1) {Slowdown (log-scale)};
    \end{tikzpicture}
    \begin{tikzpicture}[] 
    \node[] at (1,1) {\includegraphics[width=5cm]{figures/sdvbs/disparity-sqcif.pdf}};
    \node[rotate=90, align=center] at (-2, 1) {Number of Memory \\ Bandwidth Partitions};
    \node[] at (1.2, -1.8) {Number of Cache Partitions};
    \end{tikzpicture}
    \begin{tikzpicture}[]
    \node[] at (1,1) {\includegraphics[width=5cm]{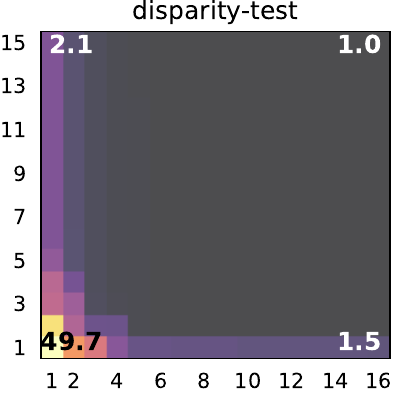}};
    \node[] at (4.5,1) {\includegraphics[height=5cm]{figures/sdvbs/disparity-sqcif_colorbar.pdf}};        
    \node[rotate=90, align=center] at (-2, 1) {Number of Memory \\ Bandwidth Partitions};
    \node[] at (1.2, -1.8) {Number of Cache Partitions};
    \node[rotate=90] at (5, 1) {Slowdown (log-scale)};
    \end{tikzpicture}
    \begin{tikzpicture}[] 
    \node[] at (1,1) {\includegraphics[width=5cm]{figures/sdvbs/disparity-vga.pdf}};
    \node[rotate=90, align=center] at (-2, 1) {Number of Memory \\ Bandwidth Partitions};
    \node[] at (1.2, -1.8) {Number of Cache Partitions};
    \end{tikzpicture}
    \begin{tikzpicture}[]
    \node[] at (1,1) {\includegraphics[width=5cm]{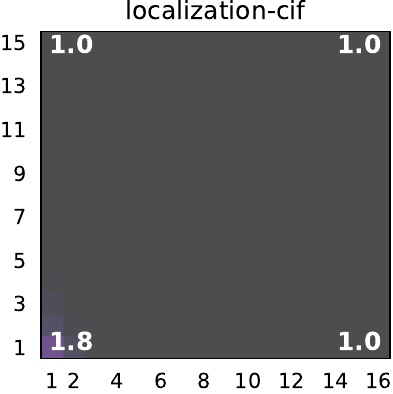}};
    \node[] at (4.5,1) {\includegraphics[height=5cm]{figures/sdvbs/disparity-sqcif_colorbar.pdf}};        
    \node[rotate=90, align=center] at (-2, 1) {Number of Memory \\ Bandwidth Partitions};
    \node[] at (1.2, -1.8) {Number of Cache Partitions};
    \node[rotate=90] at (5, 1) {Slowdown (log-scale)};
    \end{tikzpicture}
    \begin{tikzpicture}[] 
    \node[] at (1,1) {\includegraphics[width=5cm]{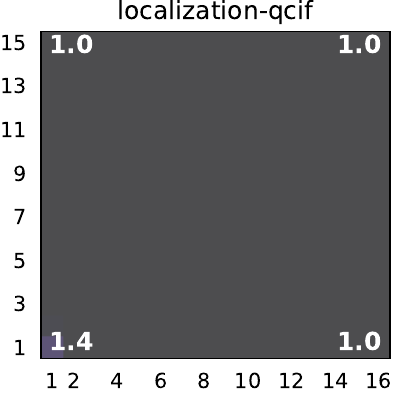}};
    \node[rotate=90, align=center] at (-2, 1) {Number of Memory \\ Bandwidth Partitions};
    \node[] at (1.2, -1.8) {Number of Cache Partitions};
    \end{tikzpicture}
    \begin{tikzpicture}[]
    \node[] at (1,1) {\includegraphics[width=5cm]{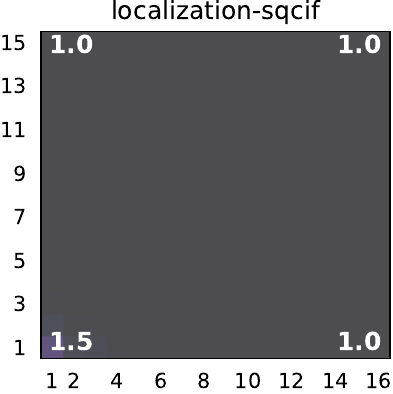}};
    \node[] at (4.5,1) {\includegraphics[height=5cm]{figures/sdvbs/disparity-sqcif_colorbar.pdf}};        
    \node[rotate=90, align=center] at (-2, 1) {Number of Memory \\ Bandwidth Partitions};
    \node[] at (1.2, -1.8) {Number of Cache Partitions};
    \node[rotate=90] at (5, 1) {Slowdown (log-scale)};
    \end{tikzpicture}
    \begin{tikzpicture}[] 
    \node[] at (1,1) {\includegraphics[width=5cm]{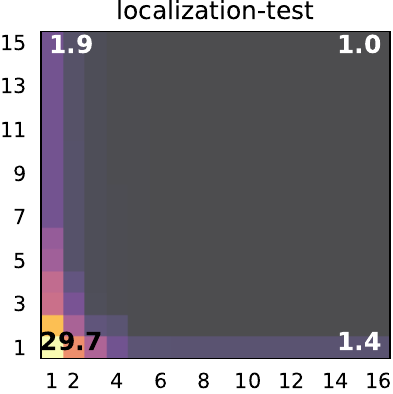}};
    \node[rotate=90, align=center] at (-2, 1) {Number of Memory \\ Bandwidth Partitions};
    \node[] at (1.2, -1.8) {Number of Cache Partitions};
    \end{tikzpicture}
    \begin{tikzpicture}[]
    \node[] at (1,1) {\includegraphics[width=5cm]{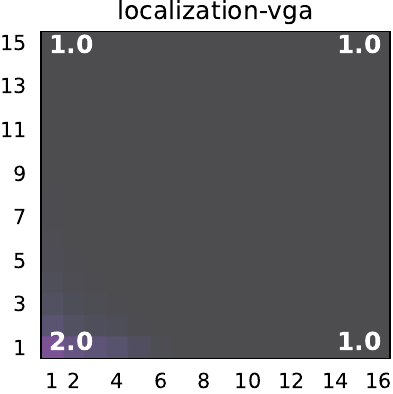}};
    \node[] at (4.5,1) {\includegraphics[height=5cm]{figures/sdvbs/disparity-sqcif_colorbar.pdf}};        
    \node[rotate=90, align=center] at (-2, 1) {Number of Memory \\ Bandwidth Partitions};
    \node[] at (1.2, -1.8) {Number of Cache Partitions};
    \node[rotate=90] at (5, 1) {Slowdown (log-scale)};
    \end{tikzpicture}
    \begin{tikzpicture}[] 
    \node[] at (1,1) {\includegraphics[width=5cm]{figures/sdvbs/mser-cif.pdf}};
    \node[rotate=90, align=center] at (-2, 1) {Number of Memory \\ Bandwidth Partitions};
    \node[] at (1.2, -1.8) {Number of Cache Partitions};
    \end{tikzpicture}
    \begin{tikzpicture}[]
    \node[] at (1,1) {\includegraphics[width=5cm]{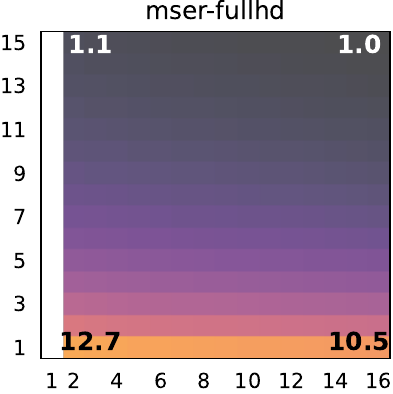}};
    \node[] at (4.5,1) {\includegraphics[height=5cm]{figures/sdvbs/disparity-sqcif_colorbar.pdf}};        
    \node[rotate=90, align=center] at (-2, 1) {Number of Memory \\ Bandwidth Partitions};
    \node[] at (1.2, -1.8) {Number of Cache Partitions};
    \node[rotate=90] at (5, 1) {Slowdown (log-scale)};
    \end{tikzpicture}
    \begin{tikzpicture}[] 
    \node[] at (1,1) {\includegraphics[width=5cm]{figures/sdvbs/mser-qcif.pdf}};
    \node[rotate=90, align=center] at (-2, 1) {Number of Memory \\ Bandwidth Partitions};
    \node[] at (1.2, -1.8) {Number of Cache Partitions};
    \end{tikzpicture}
    \begin{tikzpicture}[]
    \node[] at (1,1) {\includegraphics[width=5cm]{figures/sdvbs/mser-sqcif.pdf}};
    \node[] at (4.5,1) {\includegraphics[height=5cm]{figures/sdvbs/disparity-sqcif_colorbar.pdf}};        
    \node[rotate=90, align=center] at (-2, 1) {Number of Memory \\ Bandwidth Partitions};
    \node[] at (1.2, -1.8) {Number of Cache Partitions};
    \node[rotate=90] at (5, 1) {Slowdown (log-scale)};
    \end{tikzpicture}
    \begin{tikzpicture}[] 
    \node[] at (1,1) {\includegraphics[width=5cm]{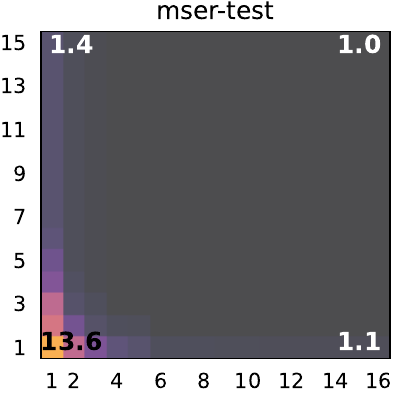}};
    \node[rotate=90, align=center] at (-2, 1) {Number of Memory \\ Bandwidth Partitions};
    \node[] at (1.2, -1.8) {Number of Cache Partitions};
    \end{tikzpicture}
    \begin{tikzpicture}[]
    \node[] at (1,1) {\includegraphics[width=5cm]{figures/sdvbs/mser-vga.pdf}};
    \node[] at (4.5,1) {\includegraphics[height=5cm]{figures/sdvbs/disparity-sqcif_colorbar.pdf}};        
    \node[rotate=90, align=center] at (-2, 1) {Number of Memory \\ Bandwidth Partitions};
    \node[] at (1.2, -1.8) {Number of Cache Partitions};
    \node[rotate=90] at (5, 1) {Slowdown (log-scale)};
    \end{tikzpicture}
    \begin{tikzpicture}[] 
    \node[] at (1,1) {\includegraphics[width=5cm]{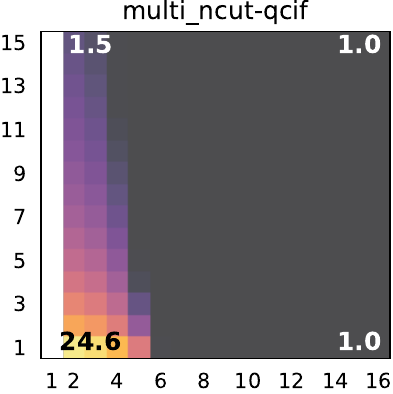}};
    \node[rotate=90, align=center] at (-2, 1) {Number of Memory \\ Bandwidth Partitions};
    \node[] at (1.2, -1.8) {Number of Cache Partitions};
    \end{tikzpicture}
    \begin{tikzpicture}[]
    \node[] at (1,1) {\includegraphics[width=5cm]{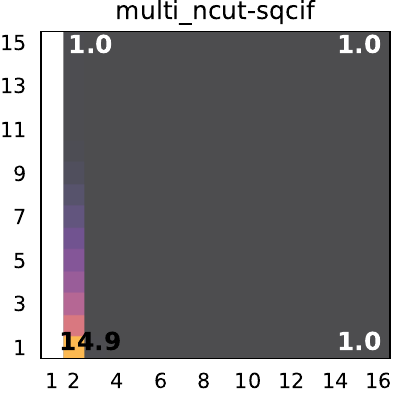}};
    \node[] at (4.5,1) {\includegraphics[height=5cm]{figures/sdvbs/disparity-sqcif_colorbar.pdf}};        
    \node[rotate=90, align=center] at (-2, 1) {Number of Memory \\ Bandwidth Partitions};
    \node[] at (1.2, -1.8) {Number of Cache Partitions};
    \node[rotate=90] at (5, 1) {Slowdown (log-scale)};
    \end{tikzpicture}
    \begin{tikzpicture}[] 
    \node[] at (1,1) {\includegraphics[width=5cm]{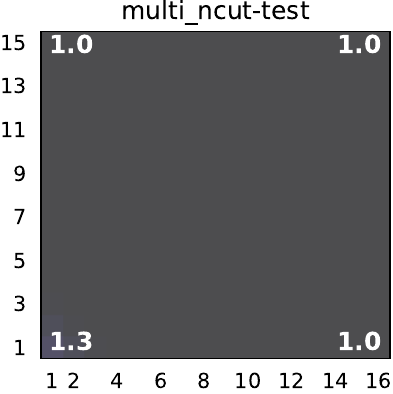}};
    \node[rotate=90, align=center] at (-2, 1) {Number of Memory \\ Bandwidth Partitions};
    \node[] at (1.2, -1.8) {Number of Cache Partitions};
    \end{tikzpicture}
    \begin{tikzpicture}[]
    \node[] at (1,1) {\includegraphics[width=5cm]{figures/sdvbs/sift-cif.pdf}};
    \node[] at (4.5,1) {\includegraphics[height=5cm]{figures/sdvbs/disparity-sqcif_colorbar.pdf}};        
    \node[rotate=90, align=center] at (-2, 1) {Number of Memory \\ Bandwidth Partitions};
    \node[] at (1.2, -1.8) {Number of Cache Partitions};
    \node[rotate=90] at (5, 1) {Slowdown (log-scale)};
    \end{tikzpicture}
    \begin{tikzpicture}[] 
    \node[] at (1,1) {\includegraphics[width=5cm]{figures/sdvbs/sift-qcif.pdf}};
    \node[rotate=90, align=center] at (-2, 1) {Number of Memory \\ Bandwidth Partitions};
    \node[] at (1.2, -1.8) {Number of Cache Partitions};
    \end{tikzpicture}
    \begin{tikzpicture}[]
    \node[] at (1,1) {\includegraphics[width=5cm]{figures/sdvbs/sift-sqcif.pdf}};
    \node[] at (4.5,1) {\includegraphics[height=5cm]{figures/sdvbs/disparity-sqcif_colorbar.pdf}};        
    \node[rotate=90, align=center] at (-2, 1) {Number of Memory \\ Bandwidth Partitions};
    \node[] at (1.2, -1.8) {Number of Cache Partitions};
    \node[rotate=90] at (5, 1) {Slowdown (log-scale)};
    \end{tikzpicture}
    \begin{tikzpicture}[] 
    \node[] at (1,1) {\includegraphics[width=5cm]{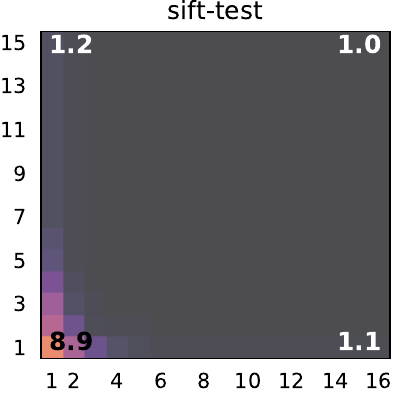}};
    \node[rotate=90, align=center] at (-2, 1) {Number of Memory \\ Bandwidth Partitions};
    \node[] at (1.2, -1.8) {Number of Cache Partitions};
    \end{tikzpicture}
    \begin{tikzpicture}[]
    \node[] at (1,1) {\includegraphics[width=5cm]{figures/sdvbs/sift-vga.pdf}};
    \node[] at (4.5,1) {\includegraphics[height=5cm]{figures/sdvbs/disparity-sqcif_colorbar.pdf}};        
    \node[rotate=90, align=center] at (-2, 1) {Number of Memory \\ Bandwidth Partitions};
    \node[] at (1.2, -1.8) {Number of Cache Partitions};
    \node[rotate=90] at (5, 1) {Slowdown (log-scale)};
    \end{tikzpicture}
    \begin{tikzpicture}[] 
    \node[] at (1,1) {\includegraphics[width=5cm]{figures/sdvbs/stitch-cif.pdf}};
    \node[rotate=90, align=center] at (-2, 1) {Number of Memory \\ Bandwidth Partitions};
    \node[] at (1.2, -1.8) {Number of Cache Partitions};
    \end{tikzpicture}
    \begin{tikzpicture}[]
    \node[] at (1,1) {\includegraphics[width=5cm]{figures/sdvbs/stitch-qcif.pdf}};
    \node[] at (4.5,1) {\includegraphics[height=5cm]{figures/sdvbs/disparity-sqcif_colorbar.pdf}};        
    \node[rotate=90, align=center] at (-2, 1) {Number of Memory \\ Bandwidth Partitions};
    \node[] at (1.2, -1.8) {Number of Cache Partitions};
    \node[rotate=90] at (5, 1) {Slowdown (log-scale)};
    \end{tikzpicture}
    \begin{tikzpicture}[] 
    \node[] at (1,1) {\includegraphics[width=5cm]{figures/sdvbs/stitch-sqcif.pdf}};
    \node[rotate=90, align=center] at (-2, 1) {Number of Memory \\ Bandwidth Partitions};
    \node[] at (1.2, -1.8) {Number of Cache Partitions};
    \end{tikzpicture}
    \begin{tikzpicture}[]
    \node[] at (1,1) {\includegraphics[width=5cm]{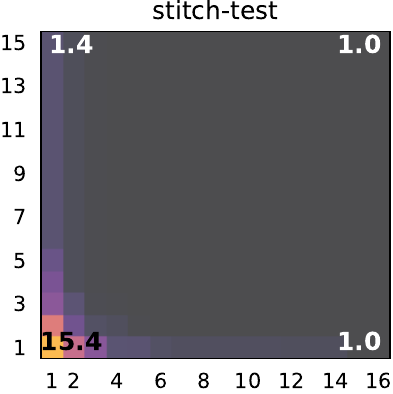}};
    \node[] at (4.5,1) {\includegraphics[height=5cm]{figures/sdvbs/disparity-sqcif_colorbar.pdf}};        
    \node[rotate=90, align=center] at (-2, 1) {Number of Memory \\ Bandwidth Partitions};
    \node[] at (1.2, -1.8) {Number of Cache Partitions};
    \node[rotate=90] at (5, 1) {Slowdown (log-scale)};
    \end{tikzpicture}
    \begin{tikzpicture}[] 
    \node[] at (1,1) {\includegraphics[width=5cm]{figures/sdvbs/stitch-vga.pdf}};
    \node[rotate=90, align=center] at (-2, 1) {Number of Memory \\ Bandwidth Partitions};
    \node[] at (1.2, -1.8) {Number of Cache Partitions};
    \end{tikzpicture}
    \begin{tikzpicture}[]
    \node[] at (1,1) {\includegraphics[width=5cm]{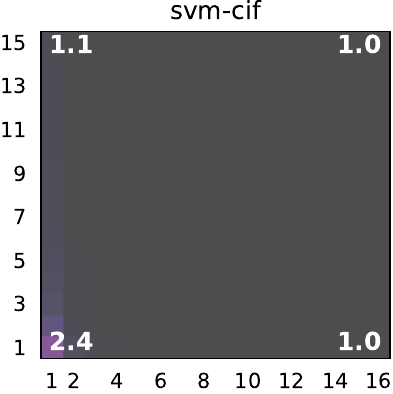}};
    \node[] at (4.5,1) {\includegraphics[height=5cm]{figures/sdvbs/disparity-sqcif_colorbar.pdf}};        
    \node[rotate=90, align=center] at (-2, 1) {Number of Memory \\ Bandwidth Partitions};
    \node[] at (1.2, -1.8) {Number of Cache Partitions};
    \node[rotate=90] at (5, 1) {Slowdown (log-scale)};
    \end{tikzpicture}
    \begin{tikzpicture}[] 
    \node[] at (1,1) {\includegraphics[width=5cm]{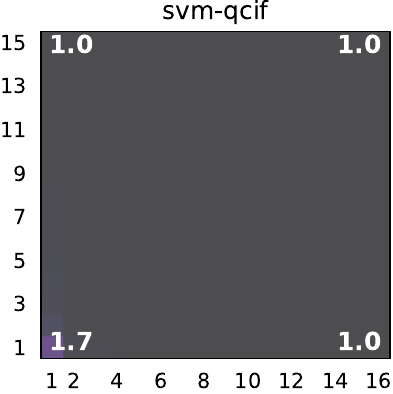}};
    \node[rotate=90, align=center] at (-2, 1) {Number of Memory \\ Bandwidth Partitions};
    \node[] at (1.2, -1.8) {Number of Cache Partitions};
    \end{tikzpicture}
    \begin{tikzpicture}[]
    \node[] at (1,1) {\includegraphics[width=5cm]{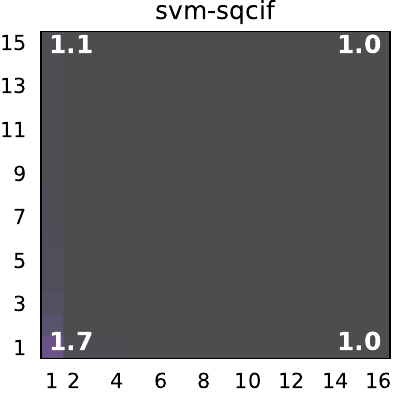}};
    \node[] at (4.5,1) {\includegraphics[height=5cm]{figures/sdvbs/disparity-sqcif_colorbar.pdf}};        
    \node[rotate=90, align=center] at (-2, 1) {Number of Memory \\ Bandwidth Partitions};
    \node[] at (1.2, -1.8) {Number of Cache Partitions};
    \node[rotate=90] at (5, 1) {Slowdown (log-scale)};
    \end{tikzpicture}
    \begin{tikzpicture}[] 
    \node[] at (1,1) {\includegraphics[width=5cm]{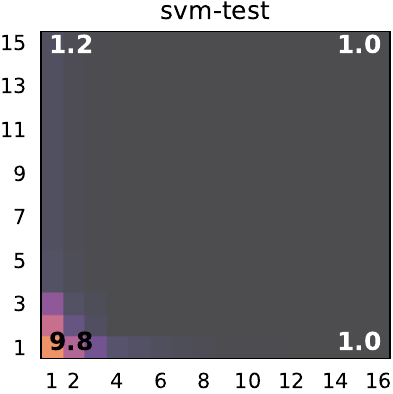}};
    \node[rotate=90, align=center] at (-2, 1) {Number of Memory \\ Bandwidth Partitions};
    \node[] at (1.2, -1.8) {Number of Cache Partitions};
    \end{tikzpicture}
    \begin{tikzpicture}[]
    \node[] at (1,1) {\includegraphics[width=5cm]{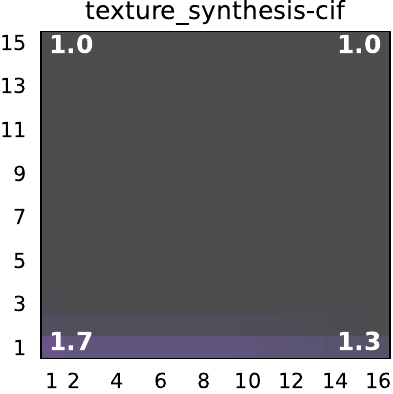}};
    \node[] at (4.5,1) {\includegraphics[height=5cm]{figures/sdvbs/disparity-sqcif_colorbar.pdf}};        
    \node[rotate=90, align=center] at (-2, 1) {Number of Memory \\ Bandwidth Partitions};
    \node[] at (1.2, -1.8) {Number of Cache Partitions};
    \node[rotate=90] at (5, 1) {Slowdown (log-scale)};
    \end{tikzpicture}
    \begin{tikzpicture}[] 
    \node[] at (1,1) {\includegraphics[width=5cm]{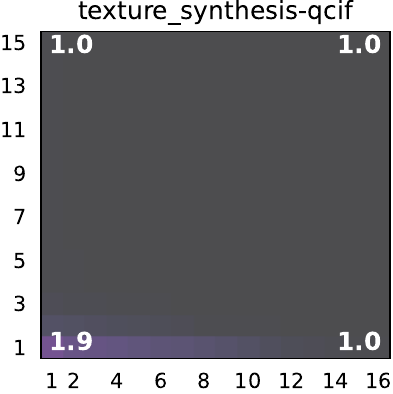}};
    \node[rotate=90, align=center] at (-2, 1) {Number of Memory \\ Bandwidth Partitions};
    \node[] at (1.2, -1.8) {Number of Cache Partitions};
    \end{tikzpicture}
    \begin{tikzpicture}[]
    \node[] at (1,1) {\includegraphics[width=5cm]{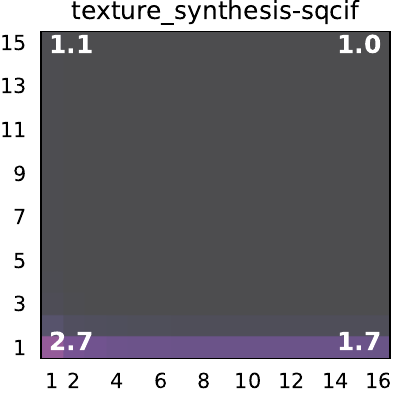}};
    \node[] at (4.5,1) {\includegraphics[height=5cm]{figures/sdvbs/disparity-sqcif_colorbar.pdf}};        
    \node[rotate=90, align=center] at (-2, 1) {Number of Memory \\ Bandwidth Partitions};
    \node[] at (1.2, -1.8) {Number of Cache Partitions};
    \node[rotate=90] at (5, 1) {Slowdown (log-scale)};
    \end{tikzpicture}
    \begin{tikzpicture}[] 
    \node[] at (1,1) {\includegraphics[width=5cm]{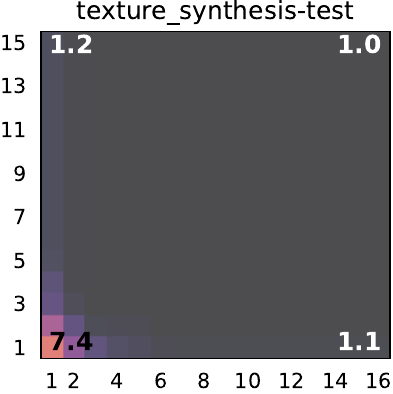}};
    \node[rotate=90, align=center] at (-2, 1) {Number of Memory \\ Bandwidth Partitions};
    \node[] at (1.2, -1.8) {Number of Cache Partitions};
    \end{tikzpicture}
    \begin{tikzpicture}[]
    \node[] at (1,1) {\includegraphics[width=5cm]{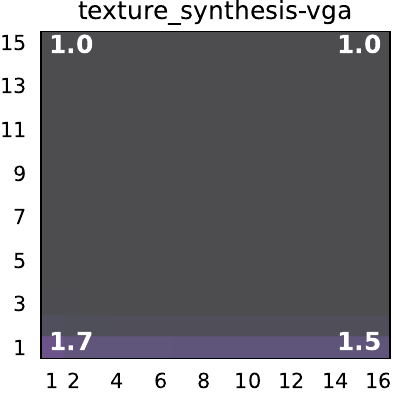}};
    \node[] at (4.5,1) {\includegraphics[height=5cm]{figures/sdvbs/disparity-sqcif_colorbar.pdf}};        
    \node[rotate=90, align=center] at (-2, 1) {Number of Memory \\ Bandwidth Partitions};
    \node[] at (1.2, -1.8) {Number of Cache Partitions};
    \node[rotate=90] at (5, 1) {Slowdown (log-scale)};
    \end{tikzpicture}
    \begin{tikzpicture}[] 
    \node[] at (1,1) {\includegraphics[width=5cm]{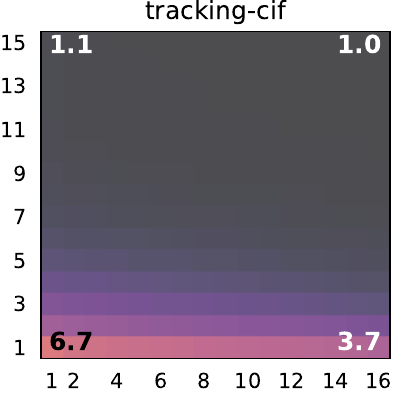}};
    \node[rotate=90, align=center] at (-2, 1) {Number of Memory \\ Bandwidth Partitions};
    \node[] at (1.2, -1.8) {Number of Cache Partitions};
    \end{tikzpicture}
    \begin{tikzpicture}[]
    \node[] at (1,1) {\includegraphics[width=5cm]{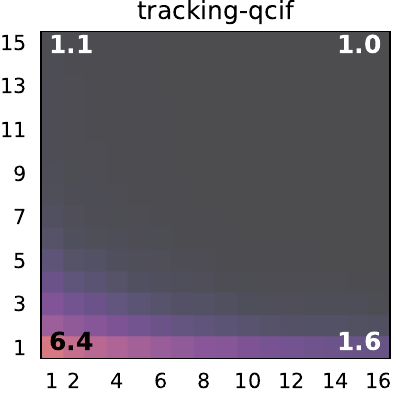}};
    \node[] at (4.5,1) {\includegraphics[height=5cm]{figures/sdvbs/disparity-sqcif_colorbar.pdf}};        
    \node[rotate=90, align=center] at (-2, 1) {Number of Memory \\ Bandwidth Partitions};
    \node[] at (1.2, -1.8) {Number of Cache Partitions};
    \node[rotate=90] at (5, 1) {Slowdown (log-scale)};
    \end{tikzpicture}
    \begin{tikzpicture}[] 
    \node[] at (1,1) {\includegraphics[width=5cm]{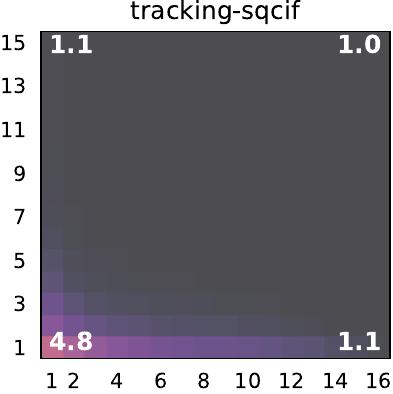}};
    \node[rotate=90, align=center] at (-2, 1) {Number of Memory \\ Bandwidth Partitions};
    \node[] at (1.2, -1.8) {Number of Cache Partitions};
    \end{tikzpicture}
    \begin{tikzpicture}[]
    \node[] at (1,1) {\includegraphics[width=5cm]{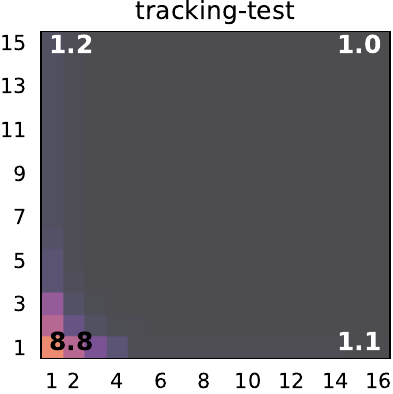}};
    \node[] at (4.5,1) {\includegraphics[height=5cm]{figures/sdvbs/disparity-sqcif_colorbar.pdf}};        
    \node[rotate=90, align=center] at (-2, 1) {Number of Memory \\ Bandwidth Partitions};
    \node[] at (1.2, -1.8) {Number of Cache Partitions};
    \node[rotate=90] at (5, 1) {Slowdown (log-scale)};
    \end{tikzpicture}
    \begin{tikzpicture}[] 
    \node[] at (1,1) {\includegraphics[width=5cm]{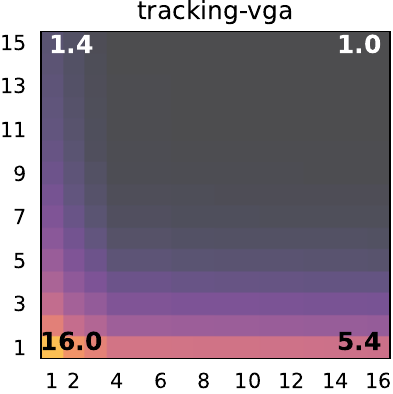}};
    \node[rotate=90, align=center] at (-2, 1) {Number of Memory \\ Bandwidth Partitions};
    \node[] at (1.2, -1.8) {Number of Cache Partitions};
    \end{tikzpicture}
    \begin{tikzpicture}[]
    \node[] at (1,1) {\includegraphics[width=5cm]{figures/sdvbs/tracking-vga.pdf}};
    \node[] at (4.5,1) {\includegraphics[height=5cm]{figures/sdvbs/disparity-sqcif_colorbar.pdf}};        
    \node[rotate=90, align=center] at (-2, 1) {Number of Memory \\ Bandwidth Partitions};
    \node[] at (1.2, -1.8) {Number of Cache Partitions};
    \node[rotate=90] at (5, 1) {Slowdown (log-scale)};
    \end{tikzpicture}

\section{PARSEC Benchmarks}

    \begin{tikzpicture}[] 
    \node[] at (1,1) {\includegraphics[width=5cm]{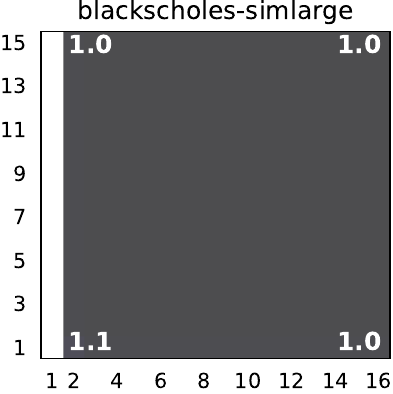}};
    \node[rotate=90, align=center] at (-2, 1) {Number of Memory \\ Bandwidth Partitions};
    \node[] at (1.2, -1.8) {Number of Cache Partitions};
    \end{tikzpicture}
    \begin{tikzpicture}[]
    \node[] at (1,1) {\includegraphics[width=5cm]{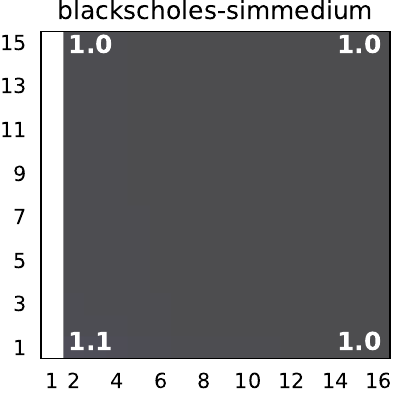}};
    \node[] at (4.5,1) {\includegraphics[height=5cm]{figures/sdvbs/disparity-sqcif_colorbar.pdf}};        
    \node[rotate=90, align=center] at (-2, 1) {Number of Memory \\ Bandwidth Partitions};
    \node[] at (1.2, -1.8) {Number of Cache Partitions};
    \node[rotate=90] at (5, 1) {Slowdown (log-scale)};
    \end{tikzpicture}
    \begin{tikzpicture}[] 
    \node[] at (1,1) {\includegraphics[width=5cm]{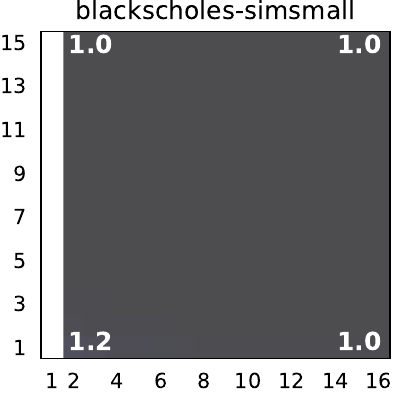}};
    \node[rotate=90, align=center] at (-2, 1) {Number of Memory \\ Bandwidth Partitions};
    \node[] at (1.2, -1.8) {Number of Cache Partitions};
    \end{tikzpicture}
    \begin{tikzpicture}[]
    \node[] at (1,1) {\includegraphics[width=5cm]{figures/parsec/bodytrack-simmedium.pdf}};
    \node[] at (4.5,1) {\includegraphics[height=5cm]{figures/sdvbs/disparity-sqcif_colorbar.pdf}};        
    \node[rotate=90, align=center] at (-2, 1) {Number of Memory \\ Bandwidth Partitions};
    \node[] at (1.2, -1.8) {Number of Cache Partitions};
    \node[rotate=90] at (5, 1) {Slowdown (log-scale)};
    \end{tikzpicture}
    \begin{tikzpicture}[] 
    \node[] at (1,1) {\includegraphics[width=5cm]{figures/parsec/bodytrack-simsmall.pdf}};
    \node[rotate=90, align=center] at (-2, 1) {Number of Memory \\ Bandwidth Partitions};
    \node[] at (1.2, -1.8) {Number of Cache Partitions};
    \end{tikzpicture}
    \begin{tikzpicture}[]
    \node[] at (1,1) {\includegraphics[width=5cm]{figures/parsec/canneal-simmedium.pdf}};
    \node[] at (4.5,1) {\includegraphics[height=5cm]{figures/sdvbs/disparity-sqcif_colorbar.pdf}};        
    \node[rotate=90, align=center] at (-2, 1) {Number of Memory \\ Bandwidth Partitions};
    \node[] at (1.2, -1.8) {Number of Cache Partitions};
    \node[rotate=90] at (5, 1) {Slowdown (log-scale)};
    \end{tikzpicture}
    \begin{tikzpicture}[] 
    \node[] at (1,1) {\includegraphics[width=5cm]{figures/parsec/canneal-simsmall.pdf}};
    \node[rotate=90, align=center] at (-2, 1) {Number of Memory \\ Bandwidth Partitions};
    \node[] at (1.2, -1.8) {Number of Cache Partitions};
    \end{tikzpicture}
    \begin{tikzpicture}[]
    \node[] at (1,1) {\includegraphics[width=5cm]{figures/parsec/ferret-simmedium.pdf}};
    \node[] at (4.5,1) {\includegraphics[height=5cm]{figures/sdvbs/disparity-sqcif_colorbar.pdf}};        
    \node[rotate=90, align=center] at (-2, 1) {Number of Memory \\ Bandwidth Partitions};
    \node[] at (1.2, -1.8) {Number of Cache Partitions};
    \node[rotate=90] at (5, 1) {Slowdown (log-scale)};
    \end{tikzpicture}
    \begin{tikzpicture}[] 
    \node[] at (1,1) {\includegraphics[width=5cm]{figures/parsec/ferret-simsmall.pdf}};
    \node[rotate=90, align=center] at (-2, 1) {Number of Memory \\ Bandwidth Partitions};
    \node[] at (1.2, -1.8) {Number of Cache Partitions};
    \end{tikzpicture}
    \begin{tikzpicture}[]
    \node[] at (1,1) {\includegraphics[width=5cm]{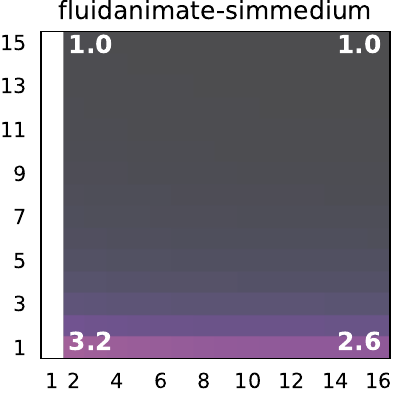}};
    \node[] at (4.5,1) {\includegraphics[height=5cm]{figures/sdvbs/disparity-sqcif_colorbar.pdf}};        
    \node[rotate=90, align=center] at (-2, 1) {Number of Memory \\ Bandwidth Partitions};
    \node[] at (1.2, -1.8) {Number of Cache Partitions};
    \node[rotate=90] at (5, 1) {Slowdown (log-scale)};
    \end{tikzpicture}
    \begin{tikzpicture}[] 
    \node[] at (1,1) {\includegraphics[width=5cm]{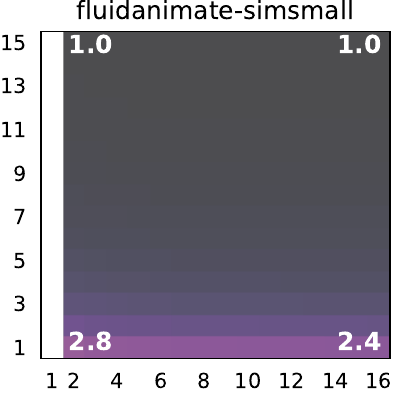}};
    \node[rotate=90, align=center] at (-2, 1) {Number of Memory \\ Bandwidth Partitions};
    \node[] at (1.2, -1.8) {Number of Cache Partitions};
    \end{tikzpicture}
    \begin{tikzpicture}[]
    \node[] at (1,1) {\includegraphics[width=5cm]{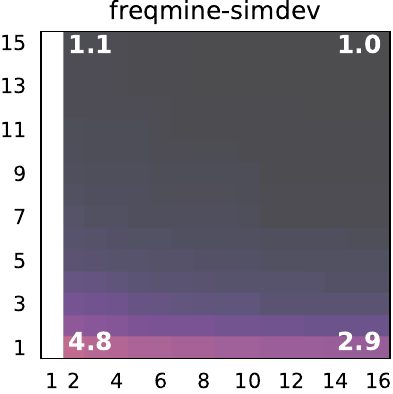}};
    \node[] at (4.5,1) {\includegraphics[height=5cm]{figures/sdvbs/disparity-sqcif_colorbar.pdf}};        
    \node[rotate=90, align=center] at (-2, 1) {Number of Memory \\ Bandwidth Partitions};
    \node[] at (1.2, -1.8) {Number of Cache Partitions};
    \node[rotate=90] at (5, 1) {Slowdown (log-scale)};
    \end{tikzpicture}
    \begin{tikzpicture}[] 
    \node[] at (1,1) {\includegraphics[width=5cm]{figures/parsec/freqmine-simmedium.pdf}};
    \node[rotate=90, align=center] at (-2, 1) {Number of Memory \\ Bandwidth Partitions};
    \node[] at (1.2, -1.8) {Number of Cache Partitions};
    \end{tikzpicture}
    \begin{tikzpicture}[]
    \node[] at (1,1) {\includegraphics[width=5cm]{figures/parsec/freqmine-simsmall.pdf}};
    \node[] at (4.5,1) {\includegraphics[height=5cm]{figures/sdvbs/disparity-sqcif_colorbar.pdf}};        
    \node[rotate=90, align=center] at (-2, 1) {Number of Memory \\ Bandwidth Partitions};
    \node[] at (1.2, -1.8) {Number of Cache Partitions};
    \node[rotate=90] at (5, 1) {Slowdown (log-scale)};
    \end{tikzpicture}
    \begin{tikzpicture}[] 
    \node[] at (1,1) {\includegraphics[width=5cm]{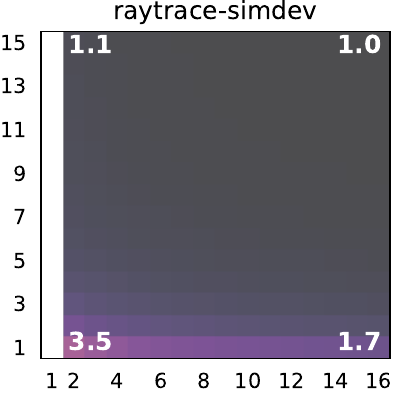}};
    \node[rotate=90, align=center] at (-2, 1) {Number of Memory \\ Bandwidth Partitions};
    \node[] at (1.2, -1.8) {Number of Cache Partitions};
    \end{tikzpicture}
    \begin{tikzpicture}[]
    \node[] at (1,1) {\includegraphics[width=5cm]{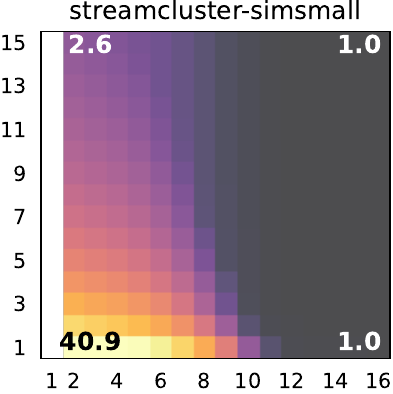}};
    \node[] at (4.5,1) {\includegraphics[height=5cm]{figures/sdvbs/disparity-sqcif_colorbar.pdf}};        
    \node[rotate=90, align=center] at (-2, 1) {Number of Memory \\ Bandwidth Partitions};
    \node[] at (1.2, -1.8) {Number of Cache Partitions};
    \node[rotate=90] at (5, 1) {Slowdown (log-scale)};
    \end{tikzpicture}
    \begin{tikzpicture}[] 
    \node[] at (1,1) {\includegraphics[width=5cm]{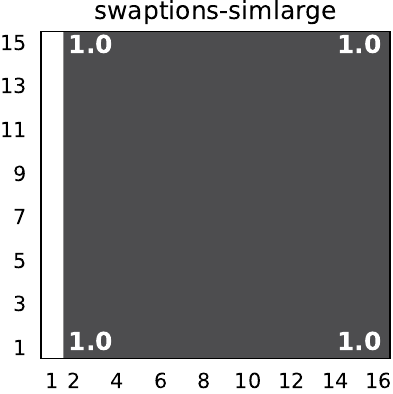}};
    \node[rotate=90, align=center] at (-2, 1) {Number of Memory \\ Bandwidth Partitions};
    \node[] at (1.2, -1.8) {Number of Cache Partitions};
    \end{tikzpicture}
    \begin{tikzpicture}[]
    \node[] at (1,1) {\includegraphics[width=5cm]{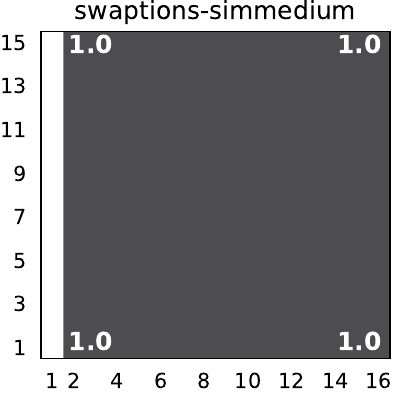}};
    \node[] at (4.5,1) {\includegraphics[height=5cm]{figures/sdvbs/disparity-sqcif_colorbar.pdf}};        
    \node[rotate=90, align=center] at (-2, 1) {Number of Memory \\ Bandwidth Partitions};
    \node[] at (1.2, -1.8) {Number of Cache Partitions};
    \node[rotate=90] at (5, 1) {Slowdown (log-scale)};
    \end{tikzpicture}
    \begin{tikzpicture}[] 
    \node[] at (1,1) {\includegraphics[width=5cm]{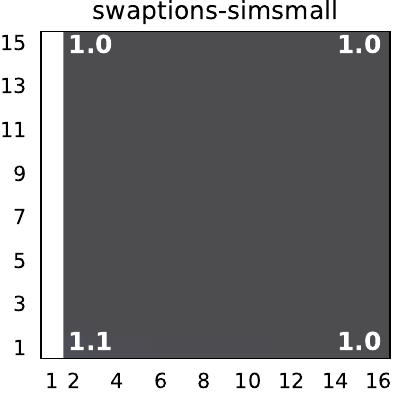}};
    \node[rotate=90, align=center] at (-2, 1) {Number of Memory \\ Bandwidth Partitions};
    \node[] at (1.2, -1.8) {Number of Cache Partitions};
    \end{tikzpicture}
    \begin{tikzpicture}[]
    \node[] at (1,1) {\includegraphics[width=5cm]{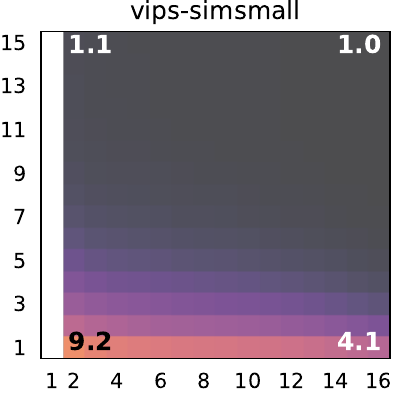}};
    \node[] at (4.5,1) {\includegraphics[height=5cm]{figures/sdvbs/disparity-sqcif_colorbar.pdf}};        
    \node[rotate=90, align=center] at (-2, 1) {Number of Memory \\ Bandwidth Partitions};
    \node[] at (1.2, -1.8) {Number of Cache Partitions};
    \node[rotate=90] at (5, 1) {Slowdown (log-scale)};
    \end{tikzpicture}

\end{appendices}
